\documentclass[11pt, twoside, reqno]{amsart}
\usepackage{amsmath,amsthm}
\usepackage{amssymb,latexsym}
\usepackage{enumerate}
\usepackage{amsfonts}
\usepackage{mathrsfs}
\usepackage{amsmath, amsthm}
\usepackage{amssymb}
\usepackage{fancyhdr}
\usepackage{esint}
\usepackage{color}
\usepackage[colorlinks,
linkcolor=red,
anchorcolor=blue,
citecolor=green]{hyperref}
\newtheorem{theorem}{\bf Theorem}[section]
\newtheorem{lemma}[theorem]{\bf Lemma}
\newtheorem{prop}[theorem]{\bf Proposition}
\newtheorem{coro}[theorem]{\bf Corollary}

\newtheorem{defn}[theorem]{\bf Definition}
\newtheorem{remark}[theorem]{\bf Remark}

\numberwithin{equation}{section}

\begin{document}
	
	\baselineskip=17pt
	
	\renewcommand{\thefootnote}{\fnsymbol {footnote}}
	
	\title[Bilinear fractional integral operators on Morrey spaces]{Bilinear fractional integral operators on Morrey spaces}
	
	\author[Qianjun He and Dunyan Yan]{Qianjun He
		\quad Dunyan Yan}
	
	\address{School of Mathematics, Graduate University, Chinese Academy of Sciences, Beijing 100049, China}
	\email{heqianjun16@mails.ucas.ac.cn}
	\email{ydunyan@ucas.ac.cn}
	\thanks{This work was supported by National Natural Science Foundation of China (Grant Nos. 11561062 and 11871452)}
	
	\begin{abstract}
	We prove a plethora of boundedness property of the Adams type for bilinear fractional integral operators of the form
	$$B_{\alpha}(f,g)(x)=\int_{\mathbb{R}^{n}}\frac{f(x-y)g(x+y)}{|y|^{n-\alpha}}dy,\qquad 0<\alpha<n.$$
	For $1<t\leq s<\infty$, we prove the non-weighted case through the known Adams type result. And we show that these results of Adams type is optimal. For $0<t\leq s<\infty$ and $0<t\leq1$, we obtain new result of a weighted theory describing Morrey boundedness of above form operators if two weights $(v,\vec{w})$ satisfy  
	\begin{equation*}
	[v,\vec{w}]_{t,\vec{q}/{a}}^{r,as}=\mathop{\sup_{Q,Q^{\prime}\in\mathscr{D}}}_{Q\subset Q^{\prime}}\left(\frac{|Q|}{|Q^{\prime}|}\right)^{\frac{1-s}{as}}|Q^{\prime}|^{\frac{1}{r}}\left(\fint_{Q}v^{\frac{t}{1-t}}\right)^{\frac{1-t}{t}}\prod_{i=1}^{2}\left(\fint_{Q^{\prime}}w_{i}^{-(q_{i}/a)^{\prime}}\right)^{\frac{1}{(q_{i}/a)^{\prime}}}<\infty,\,\,\, 0<t<s<1
	\end{equation*}
	and 
    \begin{equation*}
    [v,\vec{w}]_{t,\vec{q}/{a}}^{r,as}:=\mathop{\sup_{Q,Q^{\prime}\in\mathscr{D}}}_{Q\subset Q^{\prime}}\left(\frac{|Q|}{|Q^{\prime}|}\right)^{\frac{1-as}{as}}|Q^{\prime}|^{\frac{1}{r}}\left(\fint_{Q}v^{\frac{t}{1-t}}\right)^{\frac{1-t}{t}}\prod_{i=1}^{2}\left(\fint_{Q^{\prime}}w_{i}^{-(q_{i}/a)^{\prime}}\right)^{\frac{1}{(q_{i}/a)^{\prime}}}<\infty, \,\,\,s\geq1
    \end{equation*}
    where $\|v\|_{L^{\infty}(Q)}=\sup_{Q}v$ when $t=1$, $a$, $r$, $s$, $t$ and $\vec{q}$ satisfy proper conditions. As some applications we formulate a bilinear version of the Olsen inequality, the Fefferman-Stein type dual inequality and the Stein-Weiss inequality on Morrey spaces for fractional integrals.
	\end{abstract}

	\maketitle

	\vspace*{0.5mm}
	\section{Introduction}
	\setcounter{equation}{0}
	In the paper, we will consider the family of bilinear fractional integral operators
	\begin{equation}
	B_{\alpha}(f,g)(x):=\int_{\mathbb{R}^{n}}\frac{f(x-y)g(x+y)}{|y|^{n-\alpha}}dy,\qquad 0<\alpha<n.
	\end{equation}
	 Such operators have a long history and were studied by Bak \cite{Bak1992}, Grafakos \cite{Grafakos1992}, Grafakos and Kalton \cite{GK2001}, Hoang and Moen \cite{HM2016}, Kenig and Stein \cite{KS1999}, Kuk and Lee\cite{KL2012}, Moen \cite{Moen2014}, among others. 
	 
	For $0<\alpha<n$, the classical fractional integral $I_{\alpha}$ is given by 
	 \begin{equation}
	 I_{\alpha}f(x):=\int_{\mathbb{R}^{n}}\frac{f(y)}{|x-y|^{n-\alpha}}dy.
	 \end{equation}
	It is easily to know that $B_{\alpha}(f,g)$ and $I_{\alpha}f$ have following pointwise control relationship. For any pair of conjugate exponents ${1}/{l}+{1}/{l^{\prime}}=1$, H\"{o}lder's inequality yields
	\begin{equation}\label{Bilinear fractional integral pointwise estimate}
	|B_{\alpha}(f,g)|\lesssim I_{\alpha}(|f|^{l})^{1/l}I_{\alpha}(|g|^{l^{\prime}})^{1/l^{\prime}}.
	\end{equation}
	In \cite{Moen2014}, Moen initially introduced the fractional integral function $M_{\alpha}$, given by
	\begin{equation}\label{fractional maximal function}
	\mathcal{M}_{\alpha}(f,g)(x)=\sup_{d>0}\frac{1}{(2d)^{n-\alpha}}\int_{|y|_{\infty}\leq d}|f(x-y)g(x+y)|dy.
	\end{equation}
	A simple computation in \cite{DL2015} shows that for $0<\alpha<n$,
	$$\mathcal{M}_{\alpha}(f,g)(x)\leq cB_{\alpha}(f,g)(x).$$
	
	We first recall some stardard notation. For any measurable function $f$ the average of $f$ over a set $E$ is given by
$$\fint_{E}fdx=\frac{1}{|E|}\int_{E}fdx.$$	
The Euclidean norm of a point $x=(x_{1},\ldots,x_{n})\in\mathbb{R}^{n}$ is given by $|x|=(x_{1}^{2}+\cdots+x_{n}^{2})^{1/2}$. We also use the $l^{\infty}$ norm $|x|_{\infty}=\max(|x_{1}|,\ldots,|x_{n}|)$. Note that $|x|_{\infty}\leq|x|\leq\sqrt{n}|x|_{\infty}$ for all  $x\in\mathbb{R}^{n}$. A cube with center $x_{0}$ and side length $d$, denoted $Q=Q(x_{0},d)$, will be all points $x\in\mathbb{R}^{n}$ such that $|x-x_{0}|_{\infty}\leq\frac{d}{2}$. For an arbitrary cube $Q$, $c_{Q}$ will be its center and $l(Q)$ its side length, that is, $Q=Q(c_{Q},l(Q))$. Given $\lambda>0$ and a cube $Q$ we let $\lambda Q=Q(c_{Q},\lambda l(Q))$. The set of dyadic cubes, denoted $\mathscr{D}$, is all cubes of the form $2^{k}(m+[0,1)^{n})$ where $k\in\mathbb{Z}$ and $m\in\mathbb{Z}^{n}$. Finally for $k\in\mathbb{Z}$ we let $\mathscr{D}_{k}$ denote the cubes of level $2^{k}$, that is, $\mathscr{D}_{k}=\{Q\in\mathscr{D}:\,l(Q)=2^{k}\}$.

	Morrey spaces, named after Morrey, seem to describe the boundedness property of the classical fractional integral operators $I_{\alpha}$ more precisely than Lebesgue spaces. We first recall the definition of the Morrey (quasi-)norms \cite{Peetre1969}. For $0<q\leq p<\infty$, the Morrey norm is given by
	\begin{equation}\label{define Morrey norm}
	\|f\|_{\mathcal{M}_{q}^{p}}=\sup_{Q\in\mathscr{D}}|Q|^{{1}/{p}}\left(\fint_{Q}|f(x)|^{q}dx\right)^{\frac{1}{q}}.
	\end{equation}
	Applying H\"{o}lder's inequality to $\eqref{define Morrey norm}$, we see that
	\begin{equation}\label{the relation of Morrey norm_1}
	\|f\|_{\mathcal{M}_{q_{1}}^{p}}\geq\|f\|_{\mathcal{M}_{q_{2}}^{p}}\quad\text{for all}\,\,p\geq q_{1}\geq q_{2}>0.
	\end{equation}
	This tells us that
	\begin{equation}\label{the relation of Morrey norm_2}
	L^{p}=\mathcal{M}_{p}^{p}\subset\mathcal{M}_{q_{1}}^{p}\subset\mathcal{M}_{q_{2}}^{p}\quad\text{for all}\,\,p\geq q_{1}\geq q_{2}>0.
	\end{equation}
	\begin{remark}\label{remark Morrey spaces}
	    In addition, we know that $L^{p,\infty}$ is
		contained in $\mathcal{M}_{q}^{p}$ with $1\leq q<p<\infty$ (see \cite[Lemma 1.7]{KY1994}). More precisely, $\|f\|_{\mathcal{M}_{q}^{p}}\leq C\|f\|_{L^{p,\infty}}$ with $1\leq q<p<\infty$, here and in what follows, the letter $C$ will denote a constant, not necessarily the same in different occurrences, and let $p^{\prime}$ satisfy $1/p + 1/p^{\prime} = 1$ with $p>1$.
	\end{remark}
	
	The following result is due to Adams \cite{Adams1975} (see also Chiarenza and Frasca \cite{CF1987}), which turned out sharp \cite{Olsen1995}.
	\begin{prop}\label{Adams sharp result}
		Let $0<\alpha<n$, $1<q\leq p<\infty$ and $1<t\leq s<\infty$. Assume $\frac{1}{s}=\frac{1}{p}-\frac{\alpha}{n}$, $\frac{t}{s}=\frac{q}{p}$. Then there  exists a constant $C>0$ such that 
		$$\|I_{\alpha}f\|_{\mathcal{M}_{t}^{s}}\leq C\|f\|_{\mathcal{M}_{q}^{p}}$$
		holds for all measurable functions $f$.
	\end{prop}
		
	For the case $1<t\leq s<\infty$, we prove the following theorem under the unweighted setting.
	\begin{theorem} \label{main_1}
		Suppose that the parameters $p_{1}$, $q_{1}$, $p_{2}$, $q_{2}$, $s$, $t$ and $\alpha$ satisfy
		$$1<q_{1}\leq p_{1}<\infty,\quad1<q_{2}\leq p_{2}<\infty,\quad1/q_{1}+1/q_{2}<1,\quad1<t\leq s<\infty,\quad 0<\alpha<n.$$
		Assume that $\frac{1}{s}=\frac{1}{p_{1}}+\frac{1}{p_{2}}-\frac{\alpha}{n}$ and $\frac{t}{s}=\frac{q_{1}}{p_{1}}=\frac{q_{2}}{p_{2}}$. Then there exists a constant $C>0$ such that
		$$\|B_{\alpha}(f,g)\|_{\mathcal{M}_{t}^{s}}\leq C\|f\|_{\mathcal{M}_{q_{1}}^{p_{1}}}\|g\|_{\mathcal{M}_{q_{2}}^{p_{2}}}$$
		holds for all measurable functions $f$ and $g$.
	\end{theorem}
	
Applying inequality $\eqref{the relation of Morrey norm_1}$, we can say the following result as corollary of Theorem $\ref{main_1}$.
 
    \begin{theorem}\label{main_2}
	Suppose that the parameters $p_{1}$, $q_{1}$, $p_{2}$, $q_{2}$, $s$, $t$ and $\alpha$ satisfy 
	$$1<q_{1}\leq p_{1}<\infty,\quad1<q_{2}\leq p_{2}<\infty,\quad1/q_{1}+1/q_{2}<1,\quad1<t\leq s<\infty,\quad 0<\alpha<n.$$
	Assume that $\frac{1}{s}=\frac{1}{p_{1}}+\frac{1}{p_{2}}-\frac{\alpha}{n}$ and  $\frac{1}{t}=\frac{1}{q_{1}}+\frac{1}{q_{2}}-\frac{\alpha}{n}$. Then there exists a constant $C>0$ such that
	$$\|B_{\alpha}(f,g)\|_{\mathcal{M}_{t}^{s}}\leq C\|f\|_{\mathcal{M}_{q_{1}}^{p_{1}}}\|g\|_{\mathcal{M}_{q_{2}}^{p_{2}}}$$
	holds for all measurable functions $f$ and $g$.
    \end{theorem}
	However, this is not the end of the story; we can prove even more. Here we present our full statement of the main theorem. In specaily case Theorem $\ref*{main_1}$ can be extended to a large extent.
	\begin{theorem}\label{main_3}
		Suppose that $0<\alpha<n$, $1<q_{1}\leq p_{1}={n}/{\alpha}$, $1<q_{2}\leq p_{2}<q_{2}{n}/{\alpha}$ and ${1}/{q_{1}}+{1}/{q_{2}}<1$. Then there exists a constant $C>0$ such that
		$$\|B_{\alpha}(f,g)\|_{\mathcal{M}_{q_{2}}^{p_{2}}}\leq C\|f\|_{\mathcal{M}_{q_{1}}^{p_{1}}}\|g\|_{\mathcal{M}_{q_{2}}^{p_{2}}}$$
		holds for all positive measurable functions $f$ and $g$.
	\end{theorem}

	\section{The proofs of Theorems $\ref{main_1}-\ref{main_3}$}
	\noindent$Proof$ $of$ $Theorem$ $\ref{main_1}$. We take parameters
	$$1<u_{1},v_{1}<\infty,\qquad\qquad 1<u_{2},v_{2}<\infty,\qquad1<l<q_{1},\qquad1<l^{\prime}<q_{2}$$
	such that
	$$
	\frac{1}{u_{1}}=\frac{1}{p_{1}}-\frac{1}{l}\frac{\alpha}{n},\quad\frac{1}{u_{2}}=\frac{1}{p_{2}}-\frac{1}{l^{\prime}}\frac{\alpha}{n},\quad\frac{v_{1}}{u_{1}}=\frac{v_{2}}{u_{2}}=\frac{t}{s}.
	$$
	Since 
	$$1<l<q_{1},\quad1<l^{\prime}<q_{2}\quad\text{and}\quad\frac{1}{q_{1}}+\frac{1}{q_{2}}<1$$
	there exists a pair of conjugate of exponents $1/l+1/l^{\prime}=1$.
	
	{\noindent}Notice that 
	$$\frac{1}{u_{1}}+\frac{1}{u_{2}}=\frac{1}{p_{1}}+\frac{1}{p_{2}}-\frac{\alpha}{n}=\frac{1}{s}.$$
	It follows from this, $\frac{v_{1}}{u_{1}}=\frac{v_{2}}{u_{2}}=\frac{t}{s}$ and H\"{o}lder's inequality that
	\begin{equation}\label{Holder inequalty on Morrey spaces}
	\|h_{1}\cdot h_{2}\|_{\mathcal{M}_{t}^{s}}\leq\|h\|_{\mathcal{M}_{v_{1}}^{u_{1}}}\|h_{2}\|_{\mathcal{M}_{v_{2}}^{u_{2}}}.
	\end{equation}
	Thus, if we insert about pointwise estimate for $B_{\alpha}$ $\eqref{Bilinear fractional integral pointwise estimate}$, use the Adams original result of Proposition $\ref{Adams sharp result}$ and inequality $\eqref{Holder inequalty on Morrey spaces}$, we have 
	\begin{equation}\label{Bilinear fractional integral are controlled by two parts}
	\|B_{\alpha}(f,g)\|_{\mathcal{M}_{t}^{s}}\leq\|I_{\alpha}(|f|^{l})^{1/l}\|_{\mathcal{M}_{v_{1}}^{u_{1}}}\|I_{\alpha}(|g|^{l^{\prime}})^{1/l^{\prime}}\|_{\mathcal{M}_{v_{2}}^{u_{2}}}.
	\end{equation}
	By the condition $\frac{1}{u_{1}}=\frac{1}{p_{1}}-\frac{1}{l}\frac{\alpha}{n}$, we obtain
	\begin{equation}\label{condition equivalent_1}
	\frac{1}{u_{1}/l}=\frac{1}{p_{1}/l}-\frac{\alpha}{n}.
	\end{equation}
	Meanwhile, observing that
	\begin{equation}\label{condition equivalent_2}
		\frac{v_{1}/l}{u_{1}/l}=\frac{v_{1}}{u_{1}}=\frac{t}{s}=\frac{q_{1}}{p_{1}}=\frac{q_{1}/l}{p_{1}/l}\quad \text{and}\quad \|I_{\alpha}(|f|^{l})^{1/l}\|_{\mathcal{M}_{v_{1}}^{u_{1}}}=\|I_{\alpha}(|f|^{l})\|_{\mathcal{M}_{v_{1}/l}^{u_{1}/l}}^{1/l}.
	\end{equation}
Therefore, by equations $\eqref{condition equivalent_1}$ and $\eqref{condition equivalent_2}$ we conculde that
	\begin{equation}\label{estimate for first term}
	\|I_{\alpha}(|f|^{l})^{1/l}\|_{\mathcal{M}_{v_{1}}^{u_{1}}}\lesssim\|f^{l}\|_{\mathcal{M}_{q_{1}/l}^{p_{1}/l}}^{\frac{1}{l}}=\|f\|_{\mathcal{M}_{q_{1}}^{p_{1}}}.
	\end{equation}
	Similary, we impily that
	\begin{equation}\label{estimate for second term}
	\|I_{\alpha}(|g|^{l^{\prime}})^{1/l^{\prime}}\|_{\mathcal{M}_{v_{2}}^{u_{2}}}\lesssim \|g\|_{\mathcal{M}_{q_{2}}^{p_{2}}}.
	\end{equation}
	Combining $\eqref{Bilinear fractional integral are controlled by two parts}$, $\eqref{estimate for first term}$ and $\eqref{estimate for second term}$, we get the following estimate
	$$\|B_{\alpha}(f,g)\|_{\mathcal{M}_{t}^{s}}\leq C\|f\|_{\mathcal{M}_{q_{1}}^{p_{1}}}\|g\|_{\mathcal{M}_{q_{2}}^{p_{2}}}.$$
	This completes the proof of Theorem $\ref{main_1}$. $\hfill$ $\Box$
	
	{\noindent}$Proof$ $of$ $Theorem$ $\ref{main_2}$. Let $s$, $t_{1}$, $p_{1}$, $q_{1}$, $p_{1}$ and $q_{1}$ as in Theorem $\ref{main_1}$, then
	$$\frac{t_{1}}{s}=\frac{q_{1}}{p_{1}}=\frac{q_{2}}{p_{2}}\qquad\text{and}\qquad\frac{1}{s}=\frac{1}{p_{1}}+\frac{1}{p_{2}}-\frac{\alpha}{n}.$$
	It follows that
	$$\frac{1}{t_{1}}=\frac{p_{1}}{q_{1}}\frac{1}{s}=\frac{p_{1}}{q_{1}}\left(\frac{1}{p_{1}}+\frac{1}{p_{2}}-\frac{\alpha}{n}\right)=\frac{1}{q_{1}}+\frac{p_{1}}{q_{1}}\frac{1}{p_{2}}-\frac{p_{1}}{q_{1}}\frac{\alpha}{n}\leq\frac{1}{q_{1}}+\frac{1}{q_{1}}\frac{\alpha}{n}=\frac{1}{t},$$
	that is equivalent to 
	$$t\leq t_{1}.$$
	Therefore, by Theorem $\ref{main_1}$ and the relation $\eqref{the relation of Morrey norm_1}$ with $1\leq t\leq t_{1}<\infty$, we obtain
	$$\|B_{\alpha}(f,g)\|_{\mathcal{M}_{t}^{s}}\leq\|B_{\alpha}(f,g)\|_{\mathcal{M}_{t_{1}}^{s}}\leq C\|f\|_{\mathcal{M}_{q_{1}}^{p_{1}}}\|g\|_{\mathcal{M}_{q_{2}}^{p_{2}}}.$$
	This finishes the proof of Theorem $\ref{main_2}$. $\hfill$ $\Box$

	We invoke a bilinear estimate from \cite{SST2011}.
	\begin{prop}\label{bilinear estimate}
		Let $0<\alpha<n$, $1<p\leq p_{0}<\infty$, $1<q\leq q_{0}<\infty$ and $1<r\leq r_{0}<\infty$. Assume that
		$$q>r,\quad\frac{1}{p_{0}}>\frac{\alpha}{n},\quad\frac{1}{q_{0}}\leq\frac{\alpha}{n},$$
		and
		$$\frac{1}{r_{0}}=\frac{1}{p_{0}}+\frac{1}{q_{0}}-\frac{\alpha}{n},\quad\frac{r}{r_{0}}=\frac{p}{p_{0}}.$$
		Then
		$$\|g\cdot I_{\alpha}f\|_{\mathcal{M}_{r}^{r_{0}}}\leq C\|g\|_{\mathcal{M}_{q}^{q_{0}}}\|f\|_{\mathcal{M}_{p}^{p_{0}}},$$
		where the constant $C$ is independent of $f$ and $g$.
	\end{prop}
	
	We now prove Theorem $\ref{main_3}$.
	
	{\noindent}$Proof$ $of$ $Theorem$ $\ref{main_3}$. We first clain that we can choose parameters $1<v\leq u<\infty$ and a pair of conjugate of exponents $l,l^{\prime}>1$ such that
	\begin{equation}\label{some parameters conditions}
	p_{1}=\frac{n}{\alpha}<\frac{ln}{\alpha},\quad p_{2}<\frac{l^{\prime}n}{\alpha},\quad v>q_{2},\quad\frac{v}{u}=\frac{q_{1}}{p_{1}},\quad\frac{1}{u}=\frac{1}{p_{1}}-\frac{\alpha}{ln}=\frac{\alpha}{l^{\prime}n}.
	\end{equation}
	This is possible by assumption. In fact, let us choose  $1<v\leq u<\infty$ and $l,l^{\prime}>1$ such that
	$$
	\quad\frac{v}{u}=\frac{q_{1}}{p_{1}},\quad\frac{1}{u}=\frac{1}{p_{1}}-\frac{\alpha}{ln}.
	$$
	Then we have
	$$\frac{1}{u}=\frac{1}{p_{1}}-\frac{\alpha}{ln}=\frac{\alpha}{n}-\frac{\alpha}{ln}=\frac{\alpha}{l^{\prime}n} \quad v=\frac{q_{1}}{p_{1}}\cdot\frac{l^{\prime}n}{\alpha}.$$
	Therefore, if we choose $l,l^{\prime}$ satisfy 
	$$1<l<q_{1},\quad\max(1,\frac{q_{2}}{q_{1}})<l^{\prime}<q_{2}\quad\text{and}\quad l^{\prime}\rightarrow q_{2},$$
	Then we have
	$$v>q_{2},\quad p_{1}<\frac{ln}{\alpha}\quad\text{and}\quad p_{2}<\frac{l^{\prime}n}{\alpha}.$$
	Consequently, we could justify the claim that we can choose the parameters $1<v\leq u<\infty$ and $l,l^{\prime}>1$ so that they satisfy $\eqref{some parameters conditions}$.
	
	By inequality $\eqref{Bilinear fractional integral are controlled by two parts}$ and 
	recur to Proposition $\eqref{bilinear estimate}$ with
	$$v>q_{2},\quad p_{2}<\frac{l^{\prime}n}{\alpha},\quad u=\frac{l^{\prime}n}{\alpha},\quad \frac{1}{p_{2}}=\frac{1}{u}+\frac{1}{p_{2}}-\frac{\alpha}{l^{\prime}n}\quad\text{and}\quad\frac{q_{2}}{p_{2}}=\frac{q_{2}}{p_{2}},$$
	we have
	$$
	\|B_{\alpha}(f,g)\|_{\mathcal{M}_{q_{2}}^{p_{2}}}\leq\|I_{\alpha}(|f|^{l})^{1/l}I_{\alpha}(|g|^{l^{\prime}})^{1/l^{\prime}}\|_{\mathcal{M}_{q_{2}}^{p_{2}}}=\|I_{\alpha}(|f|^{l})^{l^{\prime}/l}\cdot I_{\alpha}(|g|^{l^{\prime}})\|_{\mathcal{M}_{q_{2}/l^{\prime}}^{p_{2}/l^{\prime}}}^{1/l^{\prime}}.
	$$
	Since
	$$\frac{v}{l^{\prime}}>\frac{q_{2}}{l^{\prime}}>1,\quad1<\frac{p_{2}}{l^{\prime}}<\frac{n}{\alpha},\quad\frac{u}{l^{\prime}}=\frac{n}{\alpha},\quad \frac{1}{p_{2}/l^{\prime}}=\frac{1}{u/l^{\prime}}+\frac{1}{p_{2}/l^{\prime}}-\frac{\alpha}{n}\quad\text{and}\quad\frac{q_{2}/l^{\prime}}{p_{2}/l^{\prime}}=\frac{q_{2}/l^{\prime}}{p_{2}/l^{\prime}},$$
	then we have
	$$
	\|B_{\alpha}(f,g)\|_{\mathcal{M}_{q_{2}}^{p_{2}}}\leq C\|I_{\alpha}(|f|^{l})^{l^{\prime}/l}\|_{\mathcal{M}_{u/l^{\prime}}^{v/l^{\prime}}}^{1/l^{\prime}}\|g\|_{\mathcal{M}_{q_{2}/l^{\prime}}^{p_{2}/l^{\prime}}}^{1/l^{\prime}}=C\|I_{\alpha}(|f|^{l})\|_{\mathcal{M}_{u/l}^{v/l}}^{1/l}\|g\|_{\mathcal{M}_{q_{2}}^{p_{2}}}.
	$$
	Meanwhile, notice that
	$$1<l<q_{1},\quad u=p_{1}l^{\prime},\quad v=q_{1}l^{\prime},\quad\frac{v/l}{u/l}=\frac{q_{1}/l}{p_{1}/l}\quad\text{and}\quad\frac{1}{u/l}=\frac{1}{p_{1}/l}-\frac{\alpha}{n}.$$
	Hence, by Proposition $\ref{Adams sharp result}$, we obtain
	$$	\|B_{\alpha}(f,g)\|_{\mathcal{M}_{q_{2}}^{p_{2}}}\leq C\|f^{l}\|_{\mathcal{M}_{q_{1}/l}^{p_{1}/l}}^{1/l}\|g\|_{\mathcal{M}_{q_{2}}^{p_{2}}}=C\|f\|_{\mathcal{M}_{q_{1}}^{p_{1}}}\|g\|_{\mathcal{M}_{q_{2}}^{p_{2}}},$$
	which gives us the desired result. $\hfill$ $\Box$

	From Theorem $\ref{main_3}$ and Remark $\ref{remark Morrey spaces}$ , we have the following result.
	\begin{coro}
		Let $\alpha$, $p_{i}$, $q_{i}$ be as in Theorem $\ref{main_3}$ and $p_{i}\neq q_{i}$ with $i=1,2$. Then
		$$\|B_{\alpha}(f,g)\|_{\mathcal{M}_{q_{2}}^{p_{2}}}\leq C\|f\|_{L^{p_{1},\infty}}\|g\|_{L^{p_{2},\infty}}.$$
	\end{coro}

	\section{Sharpness of the Results}
    In this section we prove that
    $$\|B_{\alpha}(f,g)\|_{\mathcal{M}_{t}^{s}}\leq C\|f\|_{\mathcal{M}_{q_{1}}^{p_{1}}}\|g\|_{\mathcal{M}_{q_{2}}^{p_{2}}}$$
	holds only when $\frac{t}{s}\leq\max\left(\frac{q_{1}}{p_{1}},\frac{q_{2}}{p_{2}}\right)$.
	
	\begin{theorem}\label{main_4}
		Let $0<\alpha<n$, $0<t\leq s<\infty$ and $0<q_{j}\leq p_{j}<\infty$ for $j=1,2$. Suppose that $\frac{1}{s}=\frac{1}{p_{1}}+\frac{1}{p_{2}}-\frac{\alpha}{n}$ and $\frac{t}{s}>\max\left(\frac{q_{1}}{p_{1}},\frac{q_{2}}{p_{2}}\right)$. Then there exists no constants $C>0$ such that
		$$\|B_{\alpha}(f,g)\|_{\mathcal{M}_{t}^{s}}\leq C\|f\|_{\mathcal{M}_{q_{1}}^{p_{1}}}\|g\|_{\mathcal{M}_{q_{2}}^{p_{2}}}<\infty.$$
	\end{theorem}
	$Proof$. We proof of this theorem based on following the equivalent definition of Morrey norm
\begin{equation}\label{the equivalent definition of Morrey norm}
\|f\|_{\mathcal{M}_{q}^{p}}\approx\sup_{Q\in\mathcal{Q}}|Q|^{\frac{1}{p}-\frac{1}{q}}\left(\int_{Q}|f(y)|^{q}dy\right)^{\frac{1}{q}},
\end{equation}
	where $\mathcal{Q}$ denotes the family of all open cubes in $\mathbb{R}^{n}$ with sides parallel to the coordinate axes. Without loss of generality, we may assume that
	$$1\geq\frac{t}{s}>\max\left(\frac{q_{1}}{p_{1}},\frac{q_{2}}{p_{2}}\right)=\frac{q_{1}}{p_{1}}$$
	and that $q_{1}<p_{1}$. If $1<\frac{t}{s}$, then the Morrey norm of a measurable function $f$ is infinite unless $f$ vanishes almost everywhere.
	
	Fix a positive small number $\delta\ll 1$ and $N=[\delta^{\frac{q_{1}}{p_{1}}-1}]$ be a large integer. We let the set of lattice points
	$$J:=(0,1)^{n}\cap\frac{1}{N}\mathbb{Z}^{n}.$$
	For each ponit $j\in J$, we place a small cube $Q_{j}$ centered at $j$ with the side length $\delta$ and set
	$$E:=\bigcup_{j\in J}Q_{j}.$$
	Then we have
	\begin{equation*}
	|E|=N^{n}\delta^{n}\approx\delta^{\frac{nq_{1}}{p_{1}}}.
	\end{equation*}
	Set
	$$f(x):=\delta^{-\frac{n}{p_{1}}}\chi_{3Q_{j}}(x)\quad\text{and}\quad g(x):=\delta^{-\frac{n}{p_{2}}}\chi_{3Q_{j}}(x),$$
where $3Q_{j}$ denotes the triple of $Q_{j}$.
	
	Then a simple arithmetic calculation, we claim 
\begin{equation}\label{Moreey norm of f,g }
\|f\|_{\mathcal{M}_{q_{1}}^{p_{1}}}\leq C\quad\text{and}\quad\|g\|_{\mathcal{M}_{q_{2}}^{p_{2}}}\leq C.
\end{equation}
	In fact, we use the equivalent definition of Morrey norm $\eqref{the equivalent definition of Morrey norm}$. When $l(Q)\leq3\delta$, we know that
	$$\|f\|_{\mathcal{M}_{q_{1}}^{p_{1}}}=\sup_{Q\in\mathcal{Q},l(Q)\leq3\delta}|Q|^{\frac{1}{p_{1}}-\frac{1}{q_{1}}}\delta^{-\frac{n}{p_{1}}}\left(\int_{Q\cap3Q_{j}}dy\right)^{\frac{1}{q_{1}}}\leq\sup_{Q\in\mathcal{Q},l(Q)\leq3\delta}|Q|^{\frac{1}{p_{1}}}\delta^{-\frac{n}{p_{1}}}\leq3^{\frac{n}{p_{1}}}.$$
	When $l(Q)>3\delta$, we show that
	$$\|f\|_{\mathcal{M}_{q_{1}}^{p_{1}}}=\sup_{Q\in\mathcal{Q},l(Q)>3\delta}|Q|^{\frac{1}{p_{1}}-\frac{1}{q_{1}}}\delta^{-\frac{n}{p_{1}}}\left(\int_{3Q_{j}}dy\right)^{\frac{1}{q_{1}}}\leq\sup_{Q\in\mathcal{Q},l(Q)>3\delta}3^{\frac{n}{q_{1}}}|Q|^{\frac{1}{p_{1}}-\frac{1}{q_{1}}}\delta^{\frac{1}{q_{1}}-\frac{n}{p_{1}}}\leq3^{\frac{n}{p_{1}}}.$$
	Similar to computing $\|g\|_{\mathcal{M}_{q_{2}}^{p_{2}}}$, we obtain the claim.
	
	Now, for $x\in Q_{j}$, 
	$$B_{\alpha}(f,g)(x)\approx\int_{y\in\mathbb{R}^{n}}\frac{f(x-y)g(x+y)}{|y|^{n-\alpha}}dt\geq\delta^{\alpha-n}\int_{|y|_{\infty}\leq\delta}f(x-y)g(x+y)dy=\delta^{\alpha}\delta^{-\frac{n}{p_{1}}-\frac{n}{p_{2}}}=\delta^{-\frac{n}{s}},$$
	where the above estimate based on facts: when $x\in Q$ and $|y|_{\infty}\leq l(Q)$, then $(x-y,x+y)\in 3Q\times3Q$.

	{\noindent}This tell us that
	$$\int_{(0,1)^{n}}B_{\alpha}(f,g)(x)^{t}dt\geq\int_{E}B_{\alpha}(f,g)(x)^{t}dt\leq CN^{n}\delta^{n}\delta^{-\frac{tn}{s}}\approx C(\delta^{n})^{\frac{q_{1}}{p_{1}}-\frac{t}{s}}.$$
	This implies that
	$$\|B_{\alpha}(f,g)\|_{\mathcal{M}_{t}^{s}}\geq\left(\int_{(0,1)^{n}}B_{\alpha}(f,g)(x)^{t}dt\right)^{\frac{1}{t}}\geq C\left[(\delta^{n})^{\frac{q_{1}}{p_{1}}-\frac{t}{s}}\right]^{\frac{1}{t}}.$$
	Taking $\delta$ small enough, we have the desired result by $\frac{t}{s} > \frac{q_{1}}{p_{1}}$ and $\eqref{Moreey norm of f,g }$. $\hfill$ $\Box$
	
	\section{The two-weight case for bilinear fractinal integral operators}
	
	Our results are new and provide the first non-trivial weighted estimates for $B_{\alpha}$ on Morrey spaces and the only know weighted estimates for $B_{\alpha}$ on Morrey spaces $\mathcal{M}_{t}^{s}$ when $0<t\leq1$. The estimates we obtain parallel earlier results by Iida, Sato, Sawano and Tanaka \cite{ISST2011} for the less singular bilinear fractional integral operator
	$$I_{\alpha}(f,g)(x)=\int_{\mathbb{R}^{n}}\frac{f(y)g(z)}{(|x-y|+|x-z|)^{2n-\alpha}}dydz.$$
	We first introduce two weight estimates for classical fractional integral operators on Morrey spaces, the following result is due to Iida, Sato, Sawano and Tanaka \cite{ISST2011}.
	
	\begin{prop}
		Let $v$ be a weight on $\mathbb{R}^{n}$ and $\vec{w}=(w_{1},w_{2})$ be a collection of two weights on $\mathbb{R}^{n}$. Assume that
		$$0<\alpha<2n,\quad\vec{q}=(q_{1},q_{2}),\quad 1<q_{1},q_{2}<\infty,\quad0<q\leq p<\infty,\quad0<t\leq s<r\leq\infty$$
		and $1<a<\min(r/s,q_{1},q_{2})$. Here, $q$ is given by ${1}/{q}={1}/{q_{1}}+{1}/{q_{2}}$. Suppose that
		$$\frac{1}{s}=\frac{1}{p}+\frac{1}{r}-\frac{\alpha}{n},\quad\frac{t}{s}=\frac{q}{p},\quad 0<t\leq1$$
	 and the weights ${v}$ and $\vec{w}$ satisfy the following condition:
		\begin{equation}\label{two weight condition_1}
		[v,\vec{w}]_{t,\vec{q}/{a}}^{r,as}:=\mathop{\sup_{Q,Q^{\prime}\in\mathscr{D}}}_{Q\subset Q^{\prime}}\left(\frac{|Q|}{|Q^{\prime}|}\right)^{\frac{1}{as}}|Q^{\prime}|^{\frac{1}{r}}\left(\fint_{Q}v^{t}\right)^{\frac{1}{t}}\prod_{i=1}^{2}\left(\fint_{Q^{\prime}}w_{i}^{-(q_{i}/a)^{\prime}}\right)^{\frac{1}{(q_{i}/a)^{\prime}}}<\infty.
		\end{equation}
		Then we have
		$$\|I_{\alpha}(f,g)v\|_{\mathcal{M}_{t}^{s}}\leq C	[v,\vec{w}]_{t,\vec{q}/{a}}^{r,as}\sup_{Q\in\mathscr{D}}|Q|^{1/p}\left(\fint_{Q}(|f|w_{1})^{q_{1}}\right)^{1/q_{1}}\left(\fint_{Q}(|g|w_{2})^{q_{2}}\right)^{1/q_{2}}.$$
	\end{prop}

	For the case $0<t\leq1$, we have two weight inequalities of bilinear fractional integral operators.
	\begin{theorem}\label{main_5}
			Let $v$ be a weight on $\mathbb{R}^{n}$ and $\vec{w}=(w_{1},w_{2})$ be a collection of two weights on $\mathbb{R}^{n}$. Assume that
			$$0<\alpha<n,\,\,\vec{q}=(q_{1},q_{2}),\,\, 1<q_{1},q_{2}<\infty,\,\,0<q\leq p<\infty,\,\,0<t\leq s<\infty\quad\text{and}\quad0<r\leq\infty.$$
			 Here, $q$ is given by ${1}/{q}={1}/{q_{1}}+{1}/{q_{2}}$. Suppose that
			$$\frac{\alpha}{n}>\frac{1}{r},\quad\frac{1}{s}=\frac{1}{p}+\frac{1}{r}-\frac{\alpha}{n},\quad\frac{t}{s}=\frac{q}{p},\quad0<t\leq1$$
			and the weights ${v}$ and $\vec{w}$ satisfy the following two conditions:
			
			(i) If $0<s<1$, $\frac{s}{1-s}<r$ and $1<a<\min(r(1-s)/s, q_{1}, q_{2})$
			\begin{equation}\label{two weight condition_2}
			[v,\vec{w}]_{t,\vec{q}/{a}}^{r,as}:=\mathop{\sup_{Q,Q^{\prime}\in\mathscr{D}}}_{Q\subset Q^{\prime}}\left(\frac{|Q|}{|Q^{\prime}|}\right)^{\frac{1-s}{as}}|Q^{\prime}|^{\frac{1}{r}}\left(\fint_{Q}v^{\frac{t}{1-t}}\right)^{\frac{1-t}{t}}\prod_{i=1}^{2}\left(\fint_{Q^{\prime}}w_{i}^{-(q_{i}/a)^{\prime}}\right)^{\frac{1}{(q_{i}/a)^{\prime}}}<\infty
			\end{equation}
		
		(ii) If $s\geq1$ and $1<a<\min(q_{1}, q_{2})$
			\begin{equation}\label{two weight condition_2_1}
			[v,\vec{w}]_{t,\vec{q}/{a}}^{r,as}:=\mathop{\sup_{Q,Q^{\prime}\in\mathscr{D}}}_{Q\subset Q^{\prime}}\left(\frac{|Q|}{|Q^{\prime}|}\right)^{\frac{1-as}{as}}|Q^{\prime}|^{\frac{1}{r}}\left(\fint_{Q}v^{\frac{t}{1-t}}\right)^{\frac{1-t}{t}}\prod_{i=1}^{2}\left(\fint_{Q^{\prime}}w_{i}^{-(q_{i}/a)^{\prime}}\right)^{\frac{1}{(q_{i}/a)^{\prime}}}<\infty,
			\end{equation}
			where $\left(\fint_{Q}v^{\frac{t}{1-t}}\right)^{\frac{1-t}{t}}=\|v\|_{L^{\infty}(Q)}$ when $t=1$. Then we have
			$$\|B_{\alpha}(f,g)v\|_{\mathcal{M}_{t}^{s}}\leq C	[v,\vec{w}]_{t,\vec{q}/{a}}^{r,as}\sup_{Q\in\mathscr{D}}|Q|^{1/p}\left(\fint_{Q}(|f|w_{1})^{q_{1}}\right)^{1/q_{1}}\left(\fint_{Q}(|g|w_{2})^{q_{2}}\right)^{1/q_{2}}.$$
	\end{theorem}
	
	\begin{remark}\label{remark_2}
		 Inequality $\eqref{two weight condition_2}$ holds if ${v}$ and $\vec{w}$ satisfy
		\begin{equation}\label{two condition_3}
		[v,\vec{w}]_{as,\vec{q}/{a}}^{r}=\sup_{Q\in\mathscr{D}}|Q|^{\frac{1}{r}}\left(\fint_{Q}v^{\frac{as}{1-s}}\right)^{\frac{1-s}{as}}\prod_{i=1}^{2}\left(\fint_{Q}w_{i}^{-(q_{i}/a)^{\prime}}\right)^{\frac{1}{(q_{i}/a)^{\prime}}}<\infty.
		\end{equation}
	\end{remark}
	{\noindent}Indeed, for any cubes $Q\subset Q^{\prime}$, it immediately follows that $0<t\leq s< 1$.
	Since 
	$$0<t\leq s\leq 1\Longrightarrow \frac{1-s}{as}\leq\frac{1-t}{t}\Longrightarrow\frac{t}{1-t}\leq\frac{as}{1-s},$$
	then by using H\"{o}lder's inequality we have
	\begin{align*}
	&\left(\frac{|Q|}{|Q^{\prime}|}\right)^{\frac{1-s}{as}}|Q^{\prime}|^{\frac{1}{r}}\left(\fint_{Q}v^{\frac{t}{1-t}}\right)^{\frac{1-t}{t}}\prod_{i=1}^{2}\left(\fint_{Q^{\prime}}w_{i}^{-(q_{i}/a)^{\prime}}\right)^{\frac{1}{(q_{i}/a)^{\prime}}}\\
	&\leq\left(\frac{|Q|}{|Q^{\prime}|}\right)^{\frac{1-s}{as}}|Q^{\prime}|^{\frac{1}{r}}\left(\fint_{Q}v^{\frac{as}{1-s}}\right)^{\frac{1-s}{as}}\prod_{i=1}^{2}\left(\fint_{Q^{\prime}}w_{i}^{-(q_{i}/a)^{\prime}}\right)^{\frac{1}{(q_{i}/a)^{\prime}}}\\
	&\leq |Q^{\prime}|^{\frac{1}{r}}\left(\fint_{Q^{\prime}}v^{\frac{as}{1-s}}\right)^{\frac{1-s}{as}}\prod_{i=1}^{2}\left(\fint_{Q^{\prime}}w_{i}^{-(q_{i}/a)^{\prime}}\right)^{\frac{1}{(q_{i}/a)^{\prime}}}\leq[v,\vec{w}]_{as,\vec{q}/{a}}^{r}<\infty.
	\end{align*}
	Thus, when $s=t$, $p=q$, Theorem $\ref{main_5}$ recovers the two-weight results due to Moen \cite{Moen2014}.

	The following is the Olsen inequality for bilinear fractional operators, which can see more in the papers \cite{GE2009,GSS2009,Olsen1995}.
	
	\begin{coro}
	Let $v$ be a weight on $\mathbb{R}^{n}$ and assume that
	$$0<\alpha<n,\quad 1<q_{1},q_{2}<\infty,\quad0<q\leq p<\infty,\quad0<t\leq s<1,\quad\frac{s}{1-s}<r\leq\infty\quad \text{and}\quad a>1,$$
	where $q$ is given by ${1}/{q}={1}/{q_{1}}+{1}/{q_{2}}$. Suppose that
	$$\frac{\alpha}{n}>\frac{1}{r},\quad\frac{1}{s}=\frac{1}{p}+\frac{1}{r}-\frac{\alpha}{n}\quad \text{and}\quad\frac{t}{s}=\frac{q}{p}.$$
Then we have
	$$\|B_{\alpha}(f,g)v\|_{\mathcal{M}_{t}^{s}}\leq C\|v\|_{\mathcal{M}_{\frac{t}{1-t}}^{r}}	\sup_{Q\in\mathscr{D}}|Q|^{1/p}\left(\fint_{Q}|f|^{q_{1}}\right)^{1/q_{1}}\left(\fint_{Q}|g|^{q_{2}}\right)^{1/q_{2}}.$$	
	\end{coro}
	
	$Proof$. This follows from Theorem $\ref{main_5}$ by letting $w_{1}=w_{2}=1$ and noticing that, for every $Q\subset Q^{\prime}$,
	\begin{equation}\label{set Q satisfies condition}
	\left(\frac{|Q|}{|Q^{\prime}|}\right)^{\frac{1-s}{as}}|Q^{\prime}|^{\frac{1}{r}}=|Q|^{\frac{1-s}{as}}|Q^{\prime}|^{\frac{1}{r}-\frac{1-s}{as}}\leq|Q|^{\frac{1}{r}},
	\end{equation}
	The inequality $\eqref{set Q satisfies condition}$ can be deduced from the facts that $\frac{1}{r}-\frac{1-s}{as}<0$, which follow from $\frac{s}{1-s}<r$.  
	 
	 $\hfill$ $\Box$

	The following is the Fefferman-Stein type dual inequality for bilinear fractional integrall operators on Morrey spaces.
	\begin{coro}\label{Fefferman-Stein type dual inequality}
		Assume that the parameters $0<s_{i}< 1$ and $\frac{s_{i}}{1-s_{i}}<r\leq\infty$ with $i=1,2$, satisfy
		$$
		\frac{1-s}{as}=\frac{1-s_{1}}{s_{1}}+\frac{1-s_{2}}{s_{2}}\quad\text{and}\quad\frac{1}{r}=\frac{1}{r_{1}}+\frac{1}{r_{2}}.$$
		Then, for any collection of two weights $w_{1}$ and $w_{2}$, we have
		$$\|B_{\alpha}(f,g)w_{1}w_{2}\|_{\mathcal{M}_{t}^{s}}\leq C\sup_{Q\in\mathscr{D}}|Q|^{1/p}\left(\fint_{Q}(|f|W_{1})^{q_{1}}\right)^{1/q_{1}}\left(\fint_{Q}(|g|W_{2})^{q_{2}}\right)^{1/q_{2}},
		$$
		where
		$$W_{i}(x)=\sup_{Q\in\mathscr{D}}|Q|^{1/r_{i}}\left(\fint_{Q}w_{i}^{\frac{s_{i}}{1-s_{i}}}\right)^{\frac{1-s_{i}}{s_{i}}}\quad\text{for}\,\,i=1,2.$$
	\end{coro}
	$Proof$. We need only the inequality $\eqref{two condition_3}$ with $v=w_{1}w_{2}$ and $w_{i}=W_{i}$ with $i=1,2$. It follows from H\"{o}lder's inequality that
	\begin{align*}
	Q|^{\frac{1}{r}}\left(\fint_{Q}(w_{1}w_{2})^{\frac{as}{1-s}}\right)^{\frac{1-s}{as}}\leq 	Q|^{\frac{1}{r}}\prod_{i=1}^{2}\left(\fint_{Q}w_{i}^{\frac{s_{i}}{1-s_{i}}}\right)^{\frac{1-s_{i}}{s_{i}}}=\prod_{i=1}^{2}|Q|^{1/r_{i}}\left(\fint_{Q}w_{i}^{\frac{s_{i}}{1-s_{i}}}\right)^{\frac{1-s_{i}}{s_{i}}}.
	\end{align*}
	Corollary $\ref*{Fefferman-Stein type dual inequality}$ follows immediately from the inequality
	$$W_{i}(x)\geq|Q|^{1/r_{i}}\left(\fint_{Q}w_{i}^{\frac{s_{i}}{1-s_{i}}}\right)^{\frac{1-s_{i}}{s_{i}}}\quad \text{for all}\,\,x\in Q.$$
	$\hfill$ $\Box$
	
	For one weight inequality we take $r=\infty$ and $v=w_{1}w_{2}$ to arrive at the following theorem.
	
		\begin{theorem}\label{main_6}
			Let $\vec{w}=(w_{1},w_{2})$ be a collection of two weights on $\mathbb{R}^{n}$ and assume that
			$$0<\alpha<n,\,\,\vec{q}=(q_{1},q_{2}),\,\, 1<q_{1},q_{2}<\infty,\,\,0<q\leq p<\infty,\,\,0<t\leq s<\infty \,\,\text{and}\,\, a>1,$$
			 where $q$ denotes the number determined by the H\"{o}lder relationship  ${1}/{q}={1}/{q_{1}}+{1}/{q_{2}}$. Suppose that
			$$\frac{1}{s}=\frac{1}{p}-\frac{\alpha}{n},\quad\frac{t}{s}=\frac{q}{p},\quad0<t\leq1$$
			and the weights $\vec{w}$ satisfy the following two conditions:
			
			(i) If $0<s<1$,
			\begin{equation}\label{one weight condition_1}
			[\vec{w}]_{t,\vec{q}}^{s}:=\mathop{\sup_{Q,Q^{\prime}\in\mathscr{D}}}_{Q\subset Q^{\prime}}\left(\frac{|Q|}{|Q^{\prime}|}\right)^{\frac{1-s}{as}}\left(\fint_{Q}(w_{1}w_{2})^{\frac{t}{1-t}}\right)^{\frac{1-t}{t}}\prod_{i=1}^{2}\left(\fint_{Q^{\prime}}w_{i}^{-q_{i}^{\prime}}\right)^{\frac{1}{q_{i}^{\prime}}}<\infty
			\end{equation}
		
		(ii) If $s\geq1$,
				\begin{equation}\label{one weight condition_1_1}
				[\vec{w}]_{t,\vec{q}}^{s}:=\mathop{\sup_{Q,Q^{\prime}\in\mathscr{D}}}_{Q\subset Q^{\prime}}\left(\frac{|Q|}{|Q^{\prime}|}\right)^{\frac{1-as}{as}}\left(\fint_{Q}(w_{1}w_{2})^{\frac{t}{1-t}}\right)^{\frac{1-t}{t}}\prod_{i=1}^{2}\left(\fint_{Q^{\prime}}w_{i}^{-q_{i}^{\prime}}\right)^{\frac{1}{q_{i}^{\prime}}}<\infty,
				\end{equation}
			where $\left(\fint_{Q}(w_{1}w_{2})^{\frac{t}{1-t}}\right)^{\frac{1-t}{t}}=\|w_{1}w_{2}\|_{L^{\infty}(Q)}$ when $t=1$. Then we have
			$$\|B_{\alpha}(f,g)w_{1}w_{2}\|_{\mathcal{M}_{t}^{s}}\leq C	[\vec{w}]_{t,\vec{q}}^{as}\sup_{Q\in\mathscr{D}}|Q|^{1/p}\left(\fint_{Q}(|f|w_{1})^{q_{1}}\right)^{1/q_{1}}\left(\fint_{Q}(|g|w_{2})^{q_{2}}\right)^{1/q_{2}}.$$
		\end{theorem}
		\begin{remark}\label{remark_3}
			In the same manner as in Remark $\ref{remark_2}$, by using Lemma $\ref{Characterization of a multiple weights}$ below, the inequality $\eqref{one weight condition_1}$ holds for $0<s<1$ if
		\begin{equation}\label{one weight condition_2}
		\sup_{Q\in\mathscr{D}}\left(\fint_{Q}(w_{1}w_{2})^{\frac{s}{1-s}}\right)^{\frac{1-s}{s}}\prod_{i=1}^{2}\left(\fint_{Q}w_{i}^{-q_{i}^{\prime}}\right)^{\frac{1}{q_{i}^{\prime}}}<\infty.
		\end{equation}
		\end{remark}
	{\noindent}Thus, when $s=t$ and $p=q$, Theorem $\ref{main_6}$ recovers the one-weight result due to Mone \cite{Moen2014}.

	\section{The proofs of Theorems $\ref{main_5}$ and  $\ref{main_6}$}
	We shall state and prove a principal lemma. Our key tool is the following bilinear maximal operator.
	\begin{defn}\label{define bilinear maximal operator}
		Let $0<\alpha<n$ and $0<t\leq1$. Assume that $v$ be a weight on $\mathbb{R}^{n}$ and $(f,g)$ a couple of locally integrable functions on $\mathbb{R}^{n}$. Then define a bilinear maximal operator $\widetilde{M}_{\alpha}^{t}(f,g,v)(x)$ by
		$$\widetilde{M}_{\alpha}^{t}(f,g,v)(x)=\sup_{x\in Q\in\mathscr{D}}|Q|^{\frac{\alpha}{n}}\left(\fint_{Q}|f(y)|dy\cdot\fint_{Q}|g(y)|dy\right)\left(\fint_{Q}v(y)^{\frac{t}{1-t}}dy\right)^{\frac{1-t}{t}},
		$$
		where $x\in\mathbb{R}^{n}$ and $\left(\fint_{Q}v(y)^{\frac{t}{1-t}}dy\right)^{\frac{1-t}{t}}=\|v\|_{L^{\infty}(Q)}$ when $t=1$.
	\end{defn}
	 The following is our principal lemma, which seems to be of interest on its own.
	 \begin{lemma}\label{principal lemma}
	 	Assume that $v$ be a weight on $\mathbb{R}^{n}$ and $(f,g)$ a couple of locally integrable functions on $\mathbb{R}^{n}$. For any $x\in Q_{0}\in\mathscr{D}$, set
	 	$$(f_{0}, g_{0})=(f(\cdot)\chi_{Q_{0}}(x-\cdot), g(\cdot)\chi_{Q_{0}}(\cdot-x))\quad\text{and}\quad(f_{1},g_{1})=(f\chi_{3Q_{0}},g\chi_{3Q_{0}}).$$
	 	Then there exists a constant $C$ independent of $v$, $f$, $g$ and $Q_{0}$ such that 
	 	\begin{equation}\label{bilinear fractional integral operator controlled by bilinear maximal operator }
	 	\|B_{\alpha}(f_{0},g_{0})v\|_{L^{t}(Q_{0})}\leq C\|\widetilde{M}_{\alpha}^{t}(f_{1},g_{1                                                                                                                                                                                                                 },v)\|_{L^{t}(Q_{0})}
	 	\end{equation}
	 	holds for $0<\alpha<n$ and $0<t\leq1$.
	 \end{lemma}
	 Since $B_{\alpha}$ is a positive operator, without loss of generality we may assume that $f$, $g$ are nonnegative. For simplicity, we will use the notation
	 $$m_{Q}(f,g)=\fint_{Q}f(y)dy\cdot\fint_{Q}g(y)dy.$$
	 
	 We begin with an auxilary operator that will play a key role in our analysis. For $d>0$ define,
	 $$B_{d}(f,g)(x)=\int_{|y|_{\infty}\leq d}f(x-y)g(x+y)dy.$$
	 The operators $B_{2^{k}}$ are use by Kenig and Stein \cite{KS1999} in the analysis of $B_{\alpha}$. We have the following weighted estimates for $B_{d}$ due to \cite{Moen2014}.
	 \begin{lemma}\label{Moen result}
	 		Assume that $v$ be a weight on $\mathbb{R}^{n}$ and $(f,g)$ a couple of locally integrable functions on $\mathbb{R}^{n}$. Let $0<t\leq1$ and $Q$ be a cube, then we have
	 		$$
	 		\int_{Q}B_{l(Q)}(f,g)^{t}vdx\leq C\left(\int_{3Q}fdx\cdot\int_{3Q}gdx\right)^{t}\left(\int_{Q}v^{\frac{1}{1-t}}dx\right)^{1-t},
	 		$$ 
	 		where $\left(\int_{Q}v^{\frac{1}{1-t}}dx\right)^{1-t}=\|v\|_{L^{\infty}(Q)}$ when $t=1$.
	 \end{lemma}
	 {\noindent}$Proof$. By H\"{o}lder's inequality with $1/t$ and $(1/t)^{\prime}=1/(1-t)$ we have
	 \begin{align*}
	 \int_{Q}B_{{l}(Q)}(f,g)^{t}vdx&\leq\left(\int_{Q}B_{l(Q)}(f,g)(x)dx\right)^{t}\left(\int_{Q}v^{\frac{1}{1-t}}dx\right)^{1-t}\\
	 &=\left(\int_{Q}\int_{|y|_{\infty}\leq l(Q)}f(x-y)g(x+y)dydx\right)^{t}\left(\int_{Q}v^{\frac{1}{1-t}}dx\right)^{1-t}.
	 \end{align*}
	 We make the change of variables $w=x+y$, $z=x-y$ in the first integral and notice that if $c_{Q}$ is the center of the cube, then $|x-c_{Q}|_{\infty}\leq\frac{l(Q)}{2}$ and $|t|_{\infty}\leq l(Q)$ imply that $(w,z)\in3Q\times3Q$. The lemma follows at once. $\hfill$ $\Box$
	 
	 Next we consider a discretization of the operator $B_{\alpha}$ into a dyadic model. Define the dyadic bilinear fractional integral by
	 $$B_{\alpha}^{\mathscr{D}}(f,g)(x)=\sum_{Q\in\mathscr{D}}\frac{|Q|^{\frac{\alpha}{n}}}{|Q|}B_{l(Q)}(f,g)(x)\chi_{Q}(x).$$
	
	Fix a cube $Q_{0}\in\mathscr{D}$. Let $\mathscr{D}(Q_{0})$ be the collection of all dyadic subcubes of $Q_{0}$, that is, all those cubes obtained by dividing $Q_{0}$ into $2^{n}$ congruent cubes of half its side-length, dividing each of those into $2^{n}$ congruent cubes, and so on. By convention, $Q_{0}$ itself belongs to $\mathscr{D}(Q_{0})$. To prove Lemma $\ref{principal lemma}$, we need the following estimate.
	\begin{lemma}\label{bilinear fractional integral pointwise equivalent to its dyadic model}
		For $x\in Q_{0}$, 
		\begin{equation} \label{ointwise equivalent form}
		cB_{\alpha}^{\mathscr{D}(Q_{0})}(f_{0},g_{0})(x)\leq B_{\alpha}(f_{0},g_{0})(x)\leq CB_{\alpha}^{\mathscr{D}(Q_{0})}(f_{0},g_{0})(x),
		\end{equation}
		where two constants $c$ and $C$ only depending on $\alpha$ and $n$.
	\end{lemma}
	$Proof$. We proof of $\eqref{ointwise equivalent form}$ is based on \cite{ISST2011,Moen2014}.
	
	{\noindent} We first discretize the operator $B_{\alpha}(f_{0},g_{0})$. Notice that $|y|\sim|y|_{\infty}$ and hence
	\begin{align*}
	B_{\alpha}(f_{0},g_{0})(x)&=\sum_{k\in\mathbb{Z}}\mathop{\int}_{Q_{0}\cap\{2^{k-1}<|y|_{\infty}\leq 2^{k}\}}\frac{f(x-y)g(x+y)}{|y|^{n-\alpha}}dy\\
	&\leq2^{n-\alpha}\sum_{k\in\mathbb{Z}}(2^{k})^{\alpha-n}\mathop{\int}_{Q_{0}\cap\{2^{k-1}<|y|_{\infty}\leq 2^{k}\}}\frac{f(x-y)g(x+y)}{|y|^{n-\alpha}}dy\\
	&\leq C\sum_{k\in\mathbb{Z}}(2^{k})^{\alpha-n}\mathop{\int}_{Q_{0}\cap\{|y|_{\infty}\leq 2^{k}\}}\frac{f(x-y)g(x+y)}{|y|^{n-\alpha}}dy\\
	&\leq C\sum_{x\in Q\in\mathscr{D}(Q_{0})}\frac{|Q|^{\frac{\alpha}{n}}}{|Q|}\int_{|y|_{\infty}\leq l(Q)}f(x-y)g(x+y)dy.
	\end{align*}
	
 On the other hand, fix $x\in Q_{0}$ and $\{Q_{k}\}_{k\in\mathbb{Z}}$ be the unique sequence of dyadic cubes with $x\in Q_{k}\in\mathscr{D}_{k}(Q_{0})$. Then we have
 \begin{align*}
 &\mathop{\sum_{Q\in\mathscr{D}(Q_{0})}}_{l(Q)\leq l(Q_{0})}\frac{|Q|^{\frac{\alpha}{n}}}{|Q|}B_{l(Q)}(f_{0},g_{0})(x)\chi_{Q}(x)=\sum_{k=-\infty}^{\log_{2} l(Q_{0})}\frac{|Q_{k}|^{\frac{\alpha}{n}}}{|Q_{k}|}B_{l(Q_{k})}(f_{0},g_{0})(x)\\
  &=\sum_{k=-\infty}^{\log_{2} l(Q_{0})}\frac{|Q_{k}|^{\frac{\alpha}{n}}}{|Q_{k}|}\mathop{\int}_{2^{k-1}<|y|_{\infty}\leq 2^{k}}f(x-y)g(x+y)dy+\sum_{k=-\infty}^{\log_{2} l(Q_{0})}\frac{|Q_{k}|^{\frac{\alpha}{n}}}{|Q_{k}|}B_{l(Q_{k-1})}(f_{0},g_{0})(x)\\
  &\leq c\sum_{k=-\infty}^{\log_{2} l(Q_{0})}\mathop{\int}_{2^{k-1}<|y|_{\infty}\leq 2^{k}}\frac{f(x-y)g(x+y)}{|y|^{n-\alpha}}dy+2^{\alpha-n}\mathop{\sum_{Q\in\mathscr{D}(Q_{0})}}_{l(Q)\leq l(Q_{0})}\frac{|Q|^{\frac{\alpha}{n}}}{|Q|}B_{l(Q)}(f_{0},g_{0})(x)\chi_{Q}(x)\\
  &\leq cB_{\alpha}^{\mathscr{D}(Q_{0})}(f_{0},g_{0})(x)+2^{\alpha-n}\mathop{\sum_{Q\in\mathscr{D}(Q_{0})}}_{l(Q)\leq l(Q_{0})}\frac{|Q|^{\frac{\alpha}{n}}}{|Q|}B_{l(Q)}(f_{0},g_{0})(x)\chi_{Q}(x).
 \end{align*}
Since $\alpha<n$ we may rearrange the terms, then
$$cB_{\alpha}^{\mathscr{D}(Q_{0})}(f_{0},g_{0})(x)\leq B_{\alpha}(f_{0},g_{0})(x).$$

We now proceed by following \cite{Moen2009} and observe the following. 

Define
$$M_{3\mathscr{D}}(f,g)(x)=\sup_{x\in Q\in\mathscr{D}}\fint_{3Q}fdy\cdot\fint_{3Q}gdy,$$
to be the maximal function with the basis of triples of dyadic cubes.
Letting $a>1$ be a fixed constant to be choose later, and for $k=1,2,\cdots$, we set 
$$D_{k}=\bigcup\{Q:\,Q\in\mathscr{D}(Q_{0}),\,m_{3Q}(f,g)>a^{k}\}.$$
Considering the maximal cubes with respect to inclusion, we can write
$$D_{k}=\bigcup_{j}Q_{j}^{k}$$
where the cubes $\{Q_{j}^{k}\}\subset\mathscr{D}(Q_{0})$ are nonoverlapping. By the maximality of $Q_{j}^{k}$ we can see that
\begin{equation}\label{maximality of cube}
a^{k}<m_{3Q_{j}^{k}}(f,g)\leq2^{2n}a^{k}.
\end{equation}
Let 
$$E_{0}=Q_{0}\backslash D_{1}\quad\text{and}\quad E_{j}^{k}=Q_{j}^{k}\backslash D_{k+1}.$$
We need the following properties: $\{E_{0}\}\cup\{E_{j}^{k}\}$ is a disjoint family of sets which decomposes $Q_{0}$ and satisfies
\begin{equation}\label{E_{0} and Q_{j}^{k} satisfy properties}
|Q_{0}|\leq 2|E_{0}|\quad\text{and}\quad|Q_{j}^{k}|\leq2|E_{j}^{k}|.
\end{equation}	
The inequalities $\eqref{E_{0} and Q_{j}^{k} satisfy properties}$ can be verified as follows:

Fixed $Q_{j}^{k}$ and by $\eqref{maximality of cube}$, we have that
$$Q_{j}^{k}\cap D_{k+1}\subset\{x\in Q_{j}^{k}:\,M_{3\mathscr{D}}(f,g)(x)>a^{k+1}\}.$$ 
Using the operator $M_{3\mathscr{D}}$ maps $L^{1}\times L^{1}$ into $L^{1/2,\infty}$, we have
\begin{align*}
|Q_{j}^{k}\cap D_{k+1}|&\leq|\{x\in Q_{j}^{k}:\,M_{3\mathscr{D}}(f,g)(x)>a^{k+1}\}|\\
&\leq|\{x\in\mathbb{R}^{n}:\,M_{3\mathscr{D}}(f\chi_{3Q_{j}^{k}},g\chi_{3Q_{j}^{k}})(x)>a^{k+1}\}|\\
&\leq\left(\frac{\|M_{3\mathscr{D}}\|}{a^{k+1}}\int_{3Q_{j}^{k}}f(y)dy\cdot\int_{3Q_{j}^{k}}g(y)dy\right)^{1/2}\\
&\leq\left(\frac{\|M_{3\mathscr{D}}\|}{a^{k+1}}\frac{1}{|3Q_{j}^{k}|^{2}}\int_{3Q_{j}^{k}}f(y)dy\cdot\int_{3Q_{j}^{k}}g(y)dy\right)^{1/2}|3Q_{j}^{k}|\\
&\leq\frac{6^{n}\|M_{3\mathscr{D}}\|^{1/2}}{a^{1/2}}|Q_{j}^{k}|,
\end{align*}
where $\|M_{3\mathscr{D}}\|$ be the constant from the $L^{1}\times L^{1}\rightarrow L^{1/2,\infty}$ inequality for $M_{3\mathscr{D}}$ and we have used $\eqref{maximality of cube}$ in the last step. 

Let $a=6^{2n}2^{2}\|M_{3\mathscr{D}}\|$, then we obtain
\begin{equation}\label{property condition_1}
|Q_{j}^{k}\cap D_{k+1}|\leq\frac{1}{2}|Q_{j}^{k}|.
\end{equation}
Similary, we see that 
\begin{equation}\label{property condition_2}
|D_{1}|\leq\frac{1}{2}|Q_{0}|.
\end{equation}
Clearly, $\eqref{property condition_1}$ and $\eqref{property condition_2}$ imply $\eqref{E_{0} and Q_{j}^{k} satisfy properties}$.

We set 
$$\mathscr{D}_{0}(Q_{0})=\{Q\in\mathscr{D}(Q_{0}):\,m_{3Q}(f,g)\leq a\},$$
$$\mathscr{D}_{j}^{k}(Q_{0})=\{Q\in\mathscr{D}(Q_{0}):\,Q\subset Q_{j}^{k},\,a^{k}<m_{3Q}(f,g)\leq a^{k+1}\}.$$
Then we obtain
\begin{equation}\label{Q_{0} divide into two parts}
\mathscr{D}(Q_{0})=\mathscr{D}_{0}(Q_{0})\cup\bigcup_{k,j}\mathscr{D}_{j}^{k}(Q_{0}).
\end{equation}

$Proof$ $of$ $Lemma$ $\ref{principal lemma}$. By Lemma $\ref{bilinear fractional integral pointwise equivalent to its dyadic model}$ it suffices to work the dyadic operator $B_{\alpha}^{\mathscr{D}(Q_{0})}$. Since $0<t\leq1$ it follows that
\begin{equation}\label{suffices to work the dyadic operator_1}
\int_{Q_{0}}(B_{\alpha}^{\mathscr{D}(Q_{0})}(f_{0},g_{0})v)^{t}dx\leq c\sum_{Q\in\mathscr{D}(Q_{0})}\frac{|Q|^{\frac{\alpha t}{n}}}{|Q|^{t}}\int_{Q}B_{l(Q)}(f_{0},g_{0})^{t}v^{t}dx.
\end{equation}
By Lemma $\ref{Moen result}$ we have
\begin{align}\label{suffices to work the dyadic operator_2}
&\int_{Q_{0}}(B_{\alpha}^{\mathscr{D}(Q_{0})}(f_{0},g_{0})v)^{t}dx\nonumber\\
&\leq c\sum_{Q\in\mathscr{D}(Q_{0})}\left(\frac{|Q|^{\frac{\alpha}{n}}}{|Q|}\int_{3Q}fdx\cdot\int_{3Q}gdx\right)^{t}\left(\int_{Q}v^{\frac{t}{1-t}}dx\right)^{1-t}\nonumber\\
&=c\sum_{Q\in\mathscr{D}(Q_{0})}\left({|Q|^{\frac{\alpha}{n}}}\fint_{3Q}fdx\cdot\fint_{3Q}gdx\right)^{t}\left(\fint_{Q}v^{\frac{t}{1-t}}dx\right)^{1-t}|Q|.
\end{align}
First, based on $\eqref{Q_{0} divide into two parts}$ we estimate
\begin{equation}\label{suffice to first part_1}
\sum_{Q\in\mathscr{D}_{j}^{k}(Q_{0})}\left({|Q|^{\frac{\alpha}{n}}}\fint_{3Q}fdx\cdot\fint_{3Q}gdx\right)^{t}\left(\fint_{Q}v^{\frac{t}{1-t}}dx\right)^{1-t}|Q|.
\end{equation}
For every $Q\in\mathscr{D}_{j}^{k}(Q_{0})$, we know $Q$ contained in a unique $Q_{j}^{k}$. Then
\begin{align}\label{suffice to first part_2}
\eqref{suffice to first part_1}\leq a^{(k+1)t}\sum_{Q\in\mathscr{D}_{j}^{k}(Q_{0})}|Q|^{(\frac{\alpha}{n}+1)t}\left(\int_{Q}v^{\frac{t}{1-t}}dx\right)^{1-t}\leq a^{(k+1)t}\sum_{Q\subset Q_{j}^{k}}|Q|^{(\frac{\alpha}{n}+1)t}\left(\int_{Q}v^{\frac{t}{1-t}}dx\right)^{1-t}.
\end{align}
We now use a packing condition to handle the terms in the innermost sum of $\eqref{suffice to first part_2}$. Fixe a $Q_{j}^{k}$ and consider the sum
\begin{align*}
&\sum_{Q\subset Q_{j}^{k}}|Q|^{(\frac{\alpha}{n}+1)t}\left(\int_{Q}v^{\frac{t}{1-t}}dx\right)^{1-t}\\
&=\sum_{i=0}^{\infty}\mathop{\sum_{Q\subset Q_{j}^{k}}}_{l(Q)=2^{-i}l(Q_{j}^{k})}|Q|^{(\frac{\alpha}{n}+1)t}\left(\int_{Q}v^{\frac{t}{1-t}}dx\right)^{1-t}\\
&=|Q_{j}^{k}|^{(\frac{\alpha}{n}+1)t}\sum_{i=0}^{\infty}(2^{-i\alpha t-int} )\mathop{\sum_{Q\subset Q_{j}^{k}}}_{l(Q)=2^{-i}l(Q_{j}^{k})}\left(\int_{Q}v^{\frac{t}{1-t}}dx\right)^{1-t}\\
&\leq|Q_{j}^{k}|^{(\frac{\alpha}{n}+1)t}\sum_{i=0}^{\infty}(2^{-i\alpha t-int} )\bigg(\mathop{\sum_{Q\subset Q_{j}^{k}}}_{l(Q)=2^{-i}l(Q_{j}^{k})}\int_{Q}v^{\frac{t}{1-t}}dx\bigg)^{1-t}\bigg(\mathop{\sum_{Q\subset Q_{j}^{k}}}_{l(Q)=2^{-i}l(Q_{j}^{k})}1\bigg)^{t}\\
&=|Q_{j}^{k}|^{(\frac{\alpha}{n}+1)t}\left(\int_{Q_{j}^{k}}v^{\frac{t}{1-t}}dx\right)^{1-t}\sum_{i=0}^{\infty}2^{-i\alpha t}=\frac{2^{\alpha t}}{2^{\alpha t}-1}|Q_{j}^{k}|^{\frac{\alpha}{n}t}\left(\fint_{Q_{j}^{k}}v^{\frac{t}{1-t}}dx\right)^{1-t}|Q_{j}^{k}|.
\end{align*}
Using this inequality in $\eqref{suffice to first part_2}$ we have
\begin{equation}\label{suffice to first part_3}
\eqref{suffice to first part_1}\leq C a^{(k+1)t}|Q_{j}^{k}|^{\frac{\alpha}{n}t}\left(\fint_{Q_{j}^{k}}v^{\frac{t}{1-t}}dx\right)^{1-t}|Q_{j}^{k}|.
\end{equation}
From $\eqref{maximality of cube}$, $\eqref{E_{0} and Q_{j}^{k} satisfy properties}$ and $\eqref{suffice to first part_3}$, we conclude that
\begin{align}\label{suffice to first part_4}
\eqref{suffice to first part_1}&\leq C |Q_{j}^{k}|^{\frac{\alpha}{n}t}m_{3Q_{j}^{k}}(f,g)^{t}\left(\fint_{Q_{j}^{k}}v^{\frac{t}{1-t}}dx\right)^{1-t}|Q_{j}^{k}|\nonumber\\
&=C \left[|Q_{j}^{k}|^{\frac{\alpha}{n}}m_{3Q_{j}^{k}}(f,g)\left(\fint_{Q_{j}^{k}}v^{\frac{t}{1-t}}dx\right)^{\frac{1-t}{t}}\right]^{t}|E_{j}^{k}|\leq C\int_{E_{j}^{k}}\widetilde{M}_{\alpha}^{t}(f_{0},g_{0},v)(x)^{t}dx.
\end{align}
Similary,
\begin{align}\label{suffice to first part_5}
\sum_{Q\in\mathscr{D}_{0}(Q_{0})}\left({|Q|^{\frac{\alpha}{n}}}\fint_{3Q}fdx\cdot\fint_{3Q}gdx\right)^{t}\left(\fint_{Q}v^{\frac{t}{1-t}}dx\right)^{1-t}|Q|\leq C\int_{E_{0}}\widetilde{M}_{\alpha}^{t}(f_{0},g_{0},v)(x)^{t}dx.
\end{align}
Summing up $\eqref{suffice to first part_4}$ and $\eqref{suffice to first part_5}$, we obtain
$$\int_{Q_{0}}(B_{\alpha}^{\mathscr{D}(Q_{0})}(f_{0},g_{0})v)^{t}dx\leq C\int_{Q_{0}}\widetilde{M}_{\alpha}^{t}(f_{0},g_{0},v)(x)^{t}dx.$$
This is our desired inequality $\eqref{bilinear fractional integral operator controlled by bilinear maximal operator }$. $\hfill$ $\Box$

To prove Theorems $\ref{main_5}$ and $\ref{main_6}$, we also need two more lemmas.

Let $0<\alpha<n$. For a vector $(f,g)$ of locally integrable functions and a vector $\vec{r}=(r_{1},r_{2})$ of exponents, define a maximal operator
\begin{equation}\label{define maximal operator}
M_{\alpha,\vec{r}}(f,g)(x)=\sup_{x\in Q\in\mathscr{D}}|Q|^{\frac{\alpha}{n}}\left(\fint_{Q}|f(y)|^{r_{1}}\right)^{1/r_{1}}\left(\fint_{Q}|g(y)|^{r_{2}}\right)^{1/r_{2}}.
\end{equation}

The following lemma concerns the maximal operator on Morrey spaces, which can found in the paper \cite{ISST2011}.
\begin{lemma}\label{maximal operator on Morrey spaces}
	Let $0<\alpha<n$. Set $\vec{q}=(q_{1},q_{2})$ and $\vec{r}=(r_{1},r_{2})$. Assume in addition that $0<r_{i}<q_{i}<\infty$, $i=1,2$. If $0<t\leq s<\infty$ and $0<q\leq p<\infty$ satisfy
	\begin{equation}\label{exponents condition for maximal operator on Morrey spaces}
	\frac{1}{s}=\frac{1}{p}-\frac{\alpha}{n}\quad\text{and}\quad\frac{t}{s}=\frac{q}{p},
	\end{equation}
	where $q$ is given by $1/q=1/q_{1}+1/q_{2}$, then
	$$\|M_{\alpha,\vec{r}}(f,g)\|_{\mathcal{M}_{t}^{s}}\leq C\sup_{Q\in\mathscr{D}}|Q|^{1/p}\left(\fint_{Q}|f(y)|^{q_{1}}dy\right)^{1/q_{1}}\left(\fint_{Q}|g(y)|^{q_{2}}dy\right)^{1/q_{2}}.$$
\end{lemma}

We also need the following a Characterization of a multiple weights given by Iida \cite{Iida2012}.
\begin{lemma}\label{Characterization of a multiple weights}
	Let $1<q_{1},q_{2}<\infty$ and $\hat{t}\geq q$ with $1/q=1/q_{1}+1/q_{2}$. Then, for two weights $w_{1}$, $w_{2}$, the inequality
	$$\sup_{Q\in\mathscr{D}}\left(\fint_{Q}(w_{1}w_{2})^{\hat{t}}\right)^{1/\hat{t}}\prod_{i=1}^{2}\left(\fint_{Q}w_{i}^{-q_{i}^{\prime}}\right)^{1/q_{i}^{\prime}}<\infty$$
	holds if and only if
	$$
		\left\{
		\begin{aligned}
		&(w_{1}w_{2})^{\hat{t}}\in A_{1+\hat{t}(2-1/q)}, \\
		&w_{i}^{-q_{i}^{\prime}}\in A_{q_{i}^{\prime}(1/\hat{t}+2-1/q)},\,\,i=1,2.
		\end{aligned}
		\right.
	$$
\end{lemma}

$Proof$ $of$ $Theorem$ $\ref{main_5}$. In what follows we always assume that $f$, $g$ are nonnegative and 
\begin{equation}\label{f,g of normalization}
\sup_{Q\in\mathscr{D}}|Q|^{1/p}\left(\fint_{Q}(fw_{1})^{q_{1}}\right)^{1/q_{1}}\left(\fint_{Q}(gw_{2})^{q_{2}}\right)^{1/q_{2}}=1
\end{equation}
by normalization. To prove this theorem we have estimate, for an arbitrary cube $Q_{0}\in\mathcal{Q}$, 
\begin{equation}\label{estimate for bilinear fractional integral operator on Morrey spaces}
|Q_{0}|^{1/s}\left(\fint_{Q_{0}}(B_{\alpha}(f,g)v)^{t}\right)^{1/t}.
\end{equation}
Fix a cube $Q_{0}\in\mathcal{Q}$ and recall that $(f_{0}, g_{0})=(f(\cdot)\chi_{Q_{0}}(x-\cdot), g(\cdot)\chi_{Q_{0}}(\cdot-x))$. Then by a standard argument we have, for $x\in Q_{0}$,
\begin{equation}\label{divide into two parts}
\fint_{Q_{0}}(B_{\alpha}(f,g)v)^{t}\leq\fint_{Q_{0}}(B_{0}(f_{0},g_{0})v)^{t}+C_{\infty},
\end{equation}
where 
$$C_{\infty}=\sum_{k=0}^{\infty}\fint_{Q_{0}}\left(\mathop{\int}_{2^{k}l(Q_{0})<|y|_{\infty}\leq2^{k+1}l(Q_{0})}\frac{f(x-y)g(x+y)}{|y|^{n-\alpha}}dy\right)^{t}v^{t}dx.$$
Since
$$\mathop{\int}_{2^{k}l(Q_{0})<|y|_{\infty}\leq2^{k+1}l(Q_{0})}\frac{f(x-y)g(x+y)}{|y|^{n-\alpha}}dy\leq C\frac{|2^{k+1}Q_{0}|^{\frac{\alpha}{n}}}{|2^{k+1}Q_{0}|}\int_{|y|_{\infty}\leq 2^{k+1}l(Q_{0})}f(x-y)g(x+y)dy,$$
Then we have
\begin{equation}\label{estimate for constant}
C_{\infty}\leq C\sum_{k=0}^{\infty}\frac{|2^{k+1}Q_{0}|^{\frac{\alpha t}{n}}}{|2^{k+1}Q_{0}|^{t}|Q_{0}|}\left(\int_{2^{k+3}Q_{0}}fdx\cdot\int_{2^{k+3}Q_{0}}gdx\right)^{t}\left(\int_{Q_{0}}v^{\frac{t}{1-t}}dx\right)^{1-t}.
\end{equation}

First step. Keeping in mind $\eqref{estimate for bilinear fractional integral operator on Morrey spaces}$, $\eqref{divide into two parts}$ and $\eqref{estimate for constant}$, we now estimate for $|Q_{0}|^{t/s}C_{\infty}$ in the Theorem $\ref{main_5}$. By $\eqref{two weight condition_2}$, $\eqref{two weight condition_2_1}$ and H\"{o}lder's inequality we have
$$c_{0}=\mathop{\sup_{Q\in\mathscr{D}}}_{Q_{0}\subset Q}\left(\frac{|Q_{0}|}{|Q|}\right)^{\frac{1-s}{as}}|Q|^{\frac{1}{r}}\left(\fint_{Q_{0}}v^{\frac{t}{1-t}}\right)^{\frac{1-t}{t}}\prod_{i=1}^{2}\left(\fint_{Q}w_{i}^{-q_{i}^{\prime}}\right)^{\frac{1}{q_{i}^{\prime}}}\leq[v,\vec{w}]_{t,\vec{q}/{a}}^{r,as}$$
and
$$c_{*}=\mathop{\sup_{Q\in\mathscr{D}}}_{Q_{0}\subset Q}\left(\frac{|Q_{0}|}{|Q|}\right)^{\frac{1-as}{as}}|Q|^{\frac{1}{r}}\left(\fint_{Q_{0}}v^{\frac{t}{1-t}}\right)^{\frac{1-t}{t}}\prod_{i=1}^{2}\left(\fint_{Q}w_{i}^{-q_{i}^{\prime}}\right)^{\frac{1}{q_{i}^{\prime}}}\leq[v,\vec{w}]_{t,\vec{q}/{a}}^{r,as}.$$
From H\"{o}lder's inequality, $\eqref{estimate for constant}$ and the fact that
$$\frac{1}{s}=\frac{1}{p}+\frac{1}{r}-\frac{\alpha}{n}$$
it follows that
\begin{align*}
C_{\infty}&\leq C\sum_{k=0}^{\infty}\left({\fint}_{2^{k+3}Q_{0}}(fw_{1})^{q_{1}}dx\right)^{\frac{t}{q_{1}}}\left(\fint_{2^{k+3}Q_{0}}(gw_{2})^{q_{2}}dx\right)^{\frac{t}{q_{2}}}\left(\fint_{Q_{0}}v^{\frac{t}{1-t}}dx\right)^{1-t}|2^{k+3}Q_{0}|^{\frac{\alpha t}{n}}\\
&\qquad\times\frac{|2^{k+3}Q_{0}|^{t}}{|Q_{0}|^{t}}\left(\fint_{2^{k+3}Q_{0}}w_{1}^{-q_{1}^{\prime}}\right)^{\frac{t}{q_{1}^{\prime}}}\left(\fint_{2^{k+3}Q_{0}}w_{2}^{-q_{2}^{\prime}}\right)^{\frac{t}{q_{2}^{\prime}}}\\
&=C\sum_{k=0}^{\infty}\left({\fint}_{2^{k+3}Q_{0}}(fw_{1})^{q_{1}}dx\right)^{\frac{t}{q_{1}}}\left(\fint_{2^{k+3}Q_{0}}(gw_{2})^{q_{2}}dx\right)^{\frac{t}{q_{2}}}\left(\fint_{Q_{0}}v^{\frac{t}{1-t}}dx\right)^{1-t}|2^{k+3}Q_{0}|^{\frac{t}{p}}\\
&\qquad\times\frac{|2^{k+3}Q_{0}|^{t}}{|Q_{0}|^{t}}\left(\fint_{2^{k+3}Q_{0}}w_{1}^{-q_{1}^{\prime}}\right)^{\frac{t}{q_{1}^{\prime}}}\left(\fint_{2^{k+3}Q_{0}}w_{2}^{-q_{2}^{\prime}}\right)^{\frac{t}{q_{2}^{\prime}}}|2^{k+3}Q_{0}|^{\frac{\alpha t}{n}-\frac{t}{p}}\\
&\leq C\sum_{k=0}^{\infty}\frac{|2^{k+3}Q_{0}|^{t}}{|Q_{0}|^{t}}|2^{k+3}Q_{0}|^{\frac{ t}{r}-\frac{t}{s}}\left(\fint_{2^{k+3}Q_{0}}w_{1}^{q_{1}^{\prime}}\right)^{\frac{t}{-q_{1}^{\prime}}}\left(\fint_{2^{k+3}Q_{0}}w_{2}^{-q_{2}^{\prime}}\right)^{\frac{t}{q_{2}^{\prime}}}\left(\fint_{Q_{0}}v^{\frac{t}{1-t}}dx\right)^{1-t}.
\end{align*}
This yields for $0<s<1$
$$|Q_{0}|^{t/s}C_{\infty}\leq c_{0}\sum_{k=0}^{\infty}\left(\frac{|Q_{0}|}{|2^{k+3}Q_{0}|}\right)^{t(1-1/a)(1/s-1)}=Cc_{0}$$
and for $s\geq1$
$$|Q_{0}|^{t}C_{\infty}\leq c_{0}\sum_{k=0}^{\infty}\left(\frac{|Q_{0}|}{|2^{k+3}Q_{0}|}\right)^{t(1-1/a)}=Cc_{*},$$
where we have used $1-1/a>0$.

Second step. For $0<t\leq 1$, we shall estimate
$$|Q_{0}|^{1/s}\left(\fint_{Q_{0}}(B_{\alpha}(f_{0},g_{0})v)^{t}\right)^{1/t}.$$
By $\eqref{two weight condition_2}$ and $\eqref{two weight condition_2_1}$ we have
$$c_{1}=\sup_{Q\in\mathscr{D}}|Q|^{1/r}\left(\fint_{Q}v^{\frac{t}{1-t}}\right)^{\frac{1-t}{t}}\prod_{i=1}^{2}\left(\fint_{Q}w_{i}^{-(q_{i}/a)^{\prime}}\right)^{1/(q_{i}/a)^{\prime}}\leq[v,\vec{w}]_{t,\vec{q}/a}^{r,as}.$$
To apply Lemma $\eqref{principal lemma}$ we now compute, for any $Q\in\mathscr{D}$,
\begin{align*}
&\left(\fint_{Q}|f(y)|dy\cdot\fint_{Q}|g(y)|dy\right)\left(\fint_{Q}v^{\frac{t}{1-t}}\right)^{\frac{1-t}{t}}\\
&\qquad\leq\left(\fint_{Q}(fw_{1})^{q_{1}/a}\right)^{a/q_{1}}\left(\fint_{Q}(gw_{2})^{q_{2}/a}\right)^{a/q_{2}}\left(\fint_{Q}v^{\frac{t}{1-t}}\right)^{\frac{1-t}{t}}\prod_{i=1}^{2}\left(\fint_{Q}w_{i}^{-(q_{i}/a)^{\prime}}\right)^{1/(q_{i}/a)^{\prime}}\\
&\qquad\leq c_{1}|Q|^{-1/r}\left(\fint_{Q}(fw_{1})^{q_{1}/a}\right)^{a/q_{1}}\left(\fint_{Q}(gw_{2})^{q_{2}/a}\right)^{a/q_{2}}.
\end{align*}
This implies, for $x\in Q_{0}$,
\begin{equation}\label{controlled by maximal operator}
\widetilde{M}_{\alpha}^{t}(f_{0},g_{0},v)(x)\leq c_{1}M_{\alpha-n/r,\vec{q}/a}(f,g)(x).
\end{equation}
Inequality $\eqref{controlled by maximal operator}$, Lemma $\ref{principal lemma}$ and Lemma $\ref{maximal operator on Morrey spaces}$ yield
$$|Q_{0}|^{1/s}\left(\fint_{Q_{0}}(B_{\alpha}(f,g)v)^{t}\right)^{1/t}\leq Cc_{1}\|M_{\alpha-n/r,\vec{q}/a}(f,g)\|_{\mathcal{M}_{t}^{s}}\leq Cc_{1},$$
where we have used the assumption
$$\frac{1}{s}=\frac{1}{p}-\frac{\alpha-n/r}{n}\quad\text{and}\quad\frac{t}{s}=\frac{q}{p}$$
and $\eqref{f,g of normalization}$. This completes proof of Theorem $\ref{main_5}$. $\hfill$ $\Box$

$Proof$ $of$ $Theorem$ $\ref{main_6}$.  Keeping in mind $\eqref{estimate for bilinear fractional integral operator on Morrey spaces}$, $\eqref{divide into two parts}$ and let $v=w_{1}w_{2}$, we only need estimate the following inequality
\begin{equation}\label{divide into two parts_1}
\fint_{Q_{0}}(B_{\alpha}(f,g)w_{1}w_{2})^{t}\leq\fint_{Q_{0}}(B_{0}(f_{0},g_{0})w_{1}w_{2})^{t}+c_{\infty},
\end{equation}
where 
$$c_{\infty}=\sum_{k=0}^{\infty}\fint_{Q_{0}}\left(\mathop{\int}_{2^{k}l(Q_{0})<|y|_{\infty}\leq2^{k+1}l(Q_{0})}\frac{f(x-y)g(x+y)}{|y|^{n-\alpha}}dy\right)^{t}w_{1}(x)^{t}w_{2}(x)^{t}dx.$$
Similar to the estimate for $C_{\infty}$ we have   
$$|Q_{0}|^{t/s}c_{\infty}\leq C[\vec{w}]_{t,\vec{q}}^{as}.$$
Next, we will estimate, for Theorem $\ref{main_6}$ in the conditions $\eqref{one weight condition_1}$ and $\eqref{one weight condition_1_1}$,
$$|Q_{0}|^{1/s}\left(\fint_{Q_{0}}(B_{\alpha}(f_{0},g_{0})w_{1}w_{2})^{t}\right)^{1/t}.$$
For $0<t\leq1$, by assumption we have
\begin{align*}
c_{2}=\sup_{Q\in\mathscr{D}}\left(\fint_{Q}(w_{1}w_{2})^{\frac{t}{1-t}}\right)^{\frac{1-t}{t}}\prod_{i=1}^{2}\left(\fint_{Q}w_{i}^{-q_{i}^{\prime}}\right)^{\frac{1}{q_{i}^{\prime}}}<\infty.
\end{align*}
 Then we can deduce from Lemma $\ref{Characterization of a multiple weights}$ and the reverse H\"{o}lder's inequality that there a constant $\theta\in(1,\min(q_{1},q_{2}))$  such that, for any cube $Q\in\mathscr{D}$,
\begin{equation}\label{reverse H{o}lder's inequality_1}
\left(\fint_{Q}(w_{1}w_{2})^{\frac{t}{1-t}}\right)^{\frac{1-t}{t}}\leq C\left(\fint_{Q}(w_{1}w_{2})^{\frac{t}{1-t}}\right)^{\frac{1-t}{t}}
\end{equation}
and 
\begin{equation}\label{reverse H{o}lder's inequality_2}
\left(\fint_{Q}w_{i}^{-(q_{i}/\theta)^{\prime}}\right)^{1/(q_{i}/\theta)^{\prime}}\leq C\left(\fint_{Q}w_{i}^{-q_{i}^{\prime}}\right)^{1/q_{i}^{\prime}},\quad\text{for each} \,\,i=1,2. 
\end{equation}
Combining $\eqref{reverse H{o}lder's inequality_1}$ and $\eqref{reverse H{o}lder's inequality_2}$ with the weight conditions in Theorem \ref{main_6}, we obtain 
$$[v,\vec{w}]_{t,\vec{q}/a}^{\infty,as}\leq [\vec{w}]_{t,\vec{q}}^{as}<\infty.$$
Going through a similar argument in Theorem $\ref{main_5}$, we have 
$$|Q_{0}|^{1/s}\left(\fint_{Q_{0}}(B_{\alpha}(f_{0},g_{0})w_{1}w_{2})^{t}\right)^{1/t}\leq C.$$
Consequently, Theorem $\ref{main_6}$ is proved. $\hfill$ $\Box$

\section{Examples and necessary conditions}
\noindent 7.1. {\bf{A bilinear Stein-Weiss inequality}}. Given $0<\alpha<n$ let $T_{\alpha}$ be define by 
$$T_{\alpha}f(x)=I_{n-\alpha}f(x)=\int_{\mathbb{R}^{n}}\frac{f(y)}{|x-y|^{n-\alpha}}dy.$$
Stein and Weiss\cite{SW1958} proved the following weighted inequality for $T_{\alpha}$:
$$\left(\int_{\mathbb{R}^{n}}\left(\frac{T_{\alpha}f(x)}{|x|^{\beta}}\right)^{t}dx\right)^{1/t}\leq\left(\int_{\mathbb{R}^{n}}(f(x)|x|^{\gamma})^{q}\right)^{1/q},$$
where
\begin{equation}\label{Stein-Weiss inequality condition_1}
\beta<\frac{n}{t},\qquad\gamma<\frac{n}{q^{\prime}},
\end{equation}
\begin{equation}\label{Stein-Weiss inequality condition_2}
\alpha+\beta+\gamma=n+\frac{n}{t}-\frac{n}{q},
\end{equation}
\begin{equation}\label{Stein-Weiss inequality condition_3}
\beta+\gamma\geq0.
\end{equation}
Condtions $\eqref{Stein-Weiss inequality condition_1}$, $\eqref{Stein-Weiss inequality condition_2}$ and $\eqref{Stein-Weiss inequality condition_3}$ are actually sharp. Condition $\eqref{Stein-Weiss inequality condition_1}$ ensures that $|x|^{-\beta q}$ and $|x|^{-\gamma p^{\prime}}$ are locally integrable. Condition $\eqref{Stein-Weiss inequality condition_2}$ follows from a scaling arbument and condition $\eqref{Stein-Weiss inequality condition_3}$ is a necessary condition for the weights to satisfy a general two  weight inequality \cite{SW1992}.

Below, we prove a bilinear Stein-Weiss inequality on Morrey spaces. For $0<\alpha<n$ let $BT_{\alpha}$ be the bilinear operator defined by
$$BT_{\alpha}(f,g)(x)=B_{n-\alpha}(f,g)(x)=\int_{\mathbb{R}^{n}}\frac{f(x-y)g(x+y)}{|y|^{\alpha}}dy.$$
\begin{theorem}\label{bilinear Stein-Weiss inequality on Morrey spaces}
	Assume that $1<q_{1}\leq p_{1}<\infty$, $1<q_{2}\leq p_{2}<\infty$, $0<t\leq s<1$, $\frac{n}{n-\alpha}<r\leq\infty$
	and $1<a<\min(q_{1}, q_{2})$. Here, $p$ and $q$ are given by
	$$\frac{1}{p}=\frac{1}{p_{1}}+\frac{1}{p_{2}}\quad\text{and}\quad\frac{1}{q}=\frac{1}{q_{1}}+\frac{1}{q_{2}}.$$
	Suppose that
	$$\frac{1}{s}=\frac{1}{p}+\frac{1}{r}-\frac{n-\alpha}{n},\quad\frac{1}{t}=\frac{1}{q}+\frac{1}{r}-\frac{n-\alpha}{n}$$
	 and $\alpha$, $\beta$, $\gamma_{1}$, $\gamma_{2}$ satisfy the conditons
	 \begin{equation}\label{bilinear Stein-Weiss inequality condition}
	 \left\{
	 \begin{aligned}
	 &\beta<n(\frac{1}{s}-1),\qquad \gamma_{1}<\frac{n}{q_{1}^{\prime}},\qquad \gamma_{2}<\frac{n}{q_{2}^{\prime}}, \\
	 &\quad\alpha+\beta+\gamma_{1}+\gamma_{2}=n+\frac{n}{t}-\frac{n}{q_{1}}-\frac{n}{q_{2}},\\
	 &\qquad\qquad\quad\beta+\gamma_{1}+\gamma_{2}\geq0.
	 \end{aligned}
	 \right.
	 \end{equation}
	Then the following inequality holds for all $f,g\geq 0$
	\begin{align*}
	&\sup_{Q\in\mathscr{D}}|Q|^{\frac{1}{s}}\left(\fint_{Q}\left(\frac{BT_{\alpha}(f,g)(x)}{|x|^{\beta}}\right)^{t}dx\right)^{\frac{1}{t}}\\
	&\qquad\qquad\leq C\sup_{Q\in\mathscr{D}}|Q|^{\frac{1}{p_{1}}}\left(\fint_{Q}\left(f(x)|x|^{\gamma_{1}}\right)^{q_{1}}dx\right)^{\frac{1}{q_{1}}}\sup_{Q\in\mathscr{D}}|Q|^{\frac{1}{p_{2}}}\left(\fint_{Q}\left(g(x)|x|^{\gamma_{2}}dx\right)^{q_{2}}\right)^{\frac{1}{q_{2}}}.
	\end{align*}
\end{theorem}
{\noindent}$Proof$. We suppose that the parameters satisfy
$$0<t_{0}\leq s<1\quad\text{and}\quad \frac{t_{0}}{s}=\frac{q_{1}}{p_{1}}=\frac{q_{2}}{p_{2}}.$$
Similar to the estimate for Theorem $\ref{main_2}$, we have
\begin{equation}\label{estimate for Corollary 1.4}
\sup_{Q\in\mathscr{D}}|Q|^{\frac{1}{s}}\left(\fint_{Q}\left(\frac{BT_{\alpha}(f,g)(x)}{|x|^{\beta}}\right)^{t}dx\right)^{\frac{1}{t}}\leq\sup_{Q\in\mathscr{D}}|Q|^{\frac{1}{s}}\left(\fint_{Q}\left(\frac{BT_{\alpha}(f,g)(x)}{|x|^{\beta}}\right)^{t_{0}}dx\right)^{\frac{1}{t_{0}}}.
\end{equation}
By Theorem $\ref{main_5}$, Remark $\ref{remark_2}$ and $\eqref{estimate for Corollary 1.4}$ we only need to prove that $a>1$ and a constant $C$ such that
$$|Q|^{\frac{1}{r}}\left(\fint_{Q}|x|^{-\frac{as\beta}{1-s}}dx\right)^{\frac{1-s}{as}}\prod_{i=1}^{2}\left(\fint_{Q}|x|^{-\gamma_{i}(q_{i}/a)^{\prime}}dx\right)^{\frac{1}{(q_{i}/a)^{\prime}}}\leq C$$
for all cubes $Q$. From here we follow the standard estimates for power weights. Let $a>1$ be such that $a\beta<n(\frac{1}{s}-1)$, $\gamma_{1}<\frac{n}{(q_{1}/a)^{\prime}}$ and $\gamma_{2}<\frac{n}{(q_{2}/a)^{\prime}}$. Given a cube $Q$ let $Q_{0}$ be its translate to the origin, that is, $Q_{0}=Q(0,l(Q))$. Then either $2Q_{0}\cap Q=\emptyset$ or $2Q_{0}\cap Q\neq\emptyset$. In the case $2Q_{0}\cap Q=\emptyset$ we have $|c_{Q}|_{\infty}\geq l(Q)$ and $|x|\sim|x|_{\infty}\sim|c_{Q}|_{\infty}\neq0$ for all $x\in Q$. Using this fact we have
\begin{align*}
&|Q|^{1-\frac{\alpha}{n}+\frac{1}{t}-\frac{1}{q}}\left(\fint_{Q}|x|^{-\frac{as\beta}{1-s}}dx\right)^{\frac{1-s}{as}}\prod_{i=1}^{2}\left(\fint_{Q}|x|^{-\gamma_{i}(q_{i}/a)^{\prime}}dx\right)^{\frac{1}{(q_{i}/a)^{\prime}}}\\
&=l(Q)^{\beta+\gamma_{1}+\gamma_{2}}\left(\fint_{Q}|x|^{-\frac{as\beta}{1-s}}dx\right)^{\frac{1-s}{as}}\prod_{i=1}^{2}\left(\fint_{Q}|x|^{-\gamma_{i}(q_{i}/a)^{\prime}}dx\right)^{\frac{1}{(q_{i}/a)^{\prime}}}\\
&\leq Cl(Q)^{\beta+\gamma_{1}+\gamma_{2}}|c_{Q}|_{\infty}^{-\beta-\gamma_{1}-\gamma_{2}}\leq C.
\end{align*}
where in the first line we have used the second equality in $\eqref{bilinear Stein-Weiss inequality condition}$ and in the last estimate we have used the third inequality in $\eqref{bilinear Stein-Weiss inequality condition}$. When $2Q_{0}\cap Q\neq\emptyset$ we have that $Q\subseteq B=B(0,5l(Q))$, the Euclidean ball of radius $5l(Q)$ about the origin. Thus,
\begin{align*}
&|Q|^{1-\frac{\alpha}{n}+\frac{1}{t}-\frac{1}{q}}\left(\fint_{Q}|x|^{-\frac{as\beta}{1-s}}dx\right)^{\frac{1-s}{as}}\prod_{i=1}^{2}\left(\fint_{Q}|x|^{-\gamma_{i}(q_{i}/a)^{\prime}}dx\right)^{\frac{1}{(q_{i}/a)^{\prime}}}\\
&\leq l(Q)^{\beta+\gamma_{1}+\gamma_{2}}\left(\fint_{B}|x|^{-\frac{as\beta}{1-s}}dx\right)^{\frac{1-s}{as}}\prod_{i=1}^{2}\left(\fint_{B}|x|^{-\gamma_{i}(q_{i}/a)^{\prime}}dx\right)^{\frac{1}{(q_{i}/a)^{\prime}}}\\
&\leq C.
\end{align*}
$\hfill$ $\Box$

\noindent 7.2. {\bf{Necessary conditions}}. Apparently our techniques do not address the case $1<t\leq s<\infty$. That is, other than the trivial conditions mentioned in the introduction, we do not know of sufficient conditions on weights $(v,w_{1},w_{2})$ that imply
$$
\|B_{\alpha}(f,g)v\|_{\mathcal{M}_{t}^{s}}\leq C \sup_{Q\in\mathscr{D}}|Q|^{1/p}\left(\fint_{Q}(|f|w_{1})^{q_{1}}\right)^{1/q_{1}}\left(\fint_{Q}(|g|w_{2})^{q_{2}}\right)^{1/q_{2}}
$$
when $1<t\leq s<\infty$. Here we present a necessary condition for the two weight inequality for $\mathcal{M}_{\alpha}$, which in turn is necessary for $B_{\alpha}$ when $0<\alpha<n$. 
\begin{theorem}\label{maximal operator for necessary condition}
	Let $v$ be a weight on $\mathbb{R}^{n}$ and $\vec{w}=(w_{1},w_{2})$ be a collection of two weights on $\mathbb{R}^{n}$. Assume that
	$$0\leq\alpha<n,\quad\vec{q}=(q_{1},q_{2}),\quad 1<q_{1},q_{2}<\infty,\quad0<q\leq p<\infty\quad\text{and}\quad1\leq t\leq s<\infty.$$
    Here, $q$ is given by ${1}/{q}={1}/{q_{1}}+{1}/{q_{2}}$. Suppose that
	$$\frac{\alpha}{n}\geq\frac{1}{r}\geq0,\quad\frac{1}{s}=\frac{1}{p}+\frac{1}{r}-\frac{\alpha}{n}\quad\text{and}\quad\frac{t}{s}=\frac{q}{p}.$$
    Then, for every $Q\in\mathscr{D}$, the weighted inequality
\begin{equation}\label{necessary inequality}
|Q|^{1/s}\left(\fint_{Q}(\mathcal{M}_{\alpha}(f,g)v)^{t}\right)^{1/t}\leq C\mathop{\sup_{Q^{\prime},Q\in\mathscr{D}}}_{Q\supset Q^{\prime}}|Q^{\prime}|^{1/p}\left(\fint_{Q^{\prime}}(|f|w_{1})^{q_{1}}\right)^{1/q_{1}}\left(\fint_{Q^{\prime}}(|g|w_{2})^{q_{2}}\right)^{1/q_{2}}.
\end{equation}
Then there exists a constant $C$ such that
	\begin{equation}\label{necessary condition_1}
\sup_{Q\in\mathscr{D}}|Q|^{\frac{1}{r}}(\inf_{Q}v)\left(\fint_{Q}w_{1}^{-q_{1}^{\prime}}\right)^{\frac{1}{q_{1}^{\prime}}}\left(\fint_{Q}w_{2}^{-q_{2}^{\prime}}\right)^{\frac{1}{q_{2}^{\prime}}}\leq C.
	\end{equation}
\end{theorem}

{\noindent}$Proof$. We assume to the contrary that
\begin{equation}\label{necessary condition_2}
\sup_{Q\in\mathscr{D}}|Q|^{\frac{1}{r}}(\inf_{Q}v)\left(\fint_{Q}w_{1}^{-q_{1}^{\prime}}\right)^{\frac{1}{q_{1}^{\prime}}}\left(\fint_{Q}w_{2}^{-q_{2}^{\prime}}\right)^{\frac{1}{q_{2}^{\prime}}}=\infty.
\end{equation}
By $\eqref{necessary condition_2}$ we can select a cube $Q$ such that for any large $M$,
\begin{equation}\label{necessary condition_3}
|Q|^{\frac{1}{r}}(\inf_{Q}v)\left(\fint_{Q}w_{1}^{-q_{1}^{\prime}}\right)^{\frac{1}{q_{1}^{\prime}}}\left(\fint_{Q}w_{2}^{-q_{2}^{\prime}}\right)^{\frac{1}{q_{2}^{\prime}}}>M.
\end{equation}
Selecting a smaller cube $Q$ in $\eqref{necessary condition_3}$, without loss of generality we may assume that $Q$ in minimal in the sense that
\begin{equation}\label{necessary condition_4}
\mathop{\sup_{R\in\mathscr{D}}}_{Q^{\prime}\subset R\subset Q}\fint_{R}w_{i}^{-q_{i}^{\prime}}=\fint_{Q}w_{i}^{-q_{i}^{\prime}},\quad\text{for}\,\,i=1,2.
\end{equation}
Thanks to the fact that $1/p-1/q\leq0$, equality $\eqref{necessary condition_4}$ yields
\begin{equation}\label{necessary condition_5}
\mathop{\sup_{R\in\mathscr{D}}}_{Q^{\prime}\subset R\subset Q}|R|^{1/p}\prod_{i=1}^{2}\left(\fint_{R}\chi_{Q}w_{i}^{-q_{i}^{\prime}}\right)^{1/q_{i}}=|Q|^{1/p}\prod_{i=1}^{2}\left(\fint_{Q}w_{i}^{-q_{i}^{\prime}}\right)^{1/q_{i}}.
\end{equation}
We also need the following estimate due to \cite{Moen2014},
\begin{equation}\label{Moen estimate 2014}
|Q|^{\frac{\alpha}{n}}(\sup_{Q}v)\left(\fint_{Q}f(y)dy\right)\left(\fint_{Q}g(y)dy\right)\leq C\left(\fint_{Q}(\mathcal{M}_{\alpha}(f,g)v)^{t}dx\right)^{\frac{1}{t}}.
\end{equation}
It follows by applying $\eqref{necessary inequality}$, $\eqref{necessary condition_5}$ and $\eqref{Moen estimate 2014}$ with $f=\chi_{Q}w_{1}^{-q_{1}^{\prime}}$, $g=\chi_{Q}w_{2}^{-q_{2}^{\prime}}$ that
\begin{align*}
|Q|^{\frac{1}{r}}(\inf_{Q}v)\left(\fint_{Q}w_{1}^{-q_{1}^{\prime}}\right)\left(\fint_{Q}w_{2}^{-q_{2}^{\prime}}\right)&\leq C|Q|^{-1/p}|Q|^{1/s}\left(\fint_{Q}(\mathcal{M}_{\alpha}(\chi_{Q}w_{1}^{q_{1}^{\prime}},\chi_{Q}w_{2}^{q_{2}^{\prime}})v)^{t}dx\right)^{\frac{1}{t}}\\
&\leq C|Q|^{-\frac{1}{p}}\mathop{\sup_{R\in\mathscr{D}}}_{Q\supset R\supset Q^{\prime}}|R|^{1/p}\left(\fint_{R}\chi_{Q}w_{1}^{-q_{1}^{\prime}}\right)^{1/q_{1}}\left(\fint_{R}\chi_{Q}w_{2}^{-q_{2}^{\prime}}\right)^{1/q_{2}}\\
&=C\left(\fint_{Q}w_{1}^{-q_{1}^{\prime}}\right)^{1/q_{1}}\left(\fint_{Q}w_{2}^{-q_{2}^{\prime}}\right)^{1/q_{2}}.
\end{align*}
This yields a contradiction
$$M<|Q|^{\frac{1}{r}}(\inf_{Q}v)\left(\fint_{Q}w_{1}^{-q_{1}^{\prime}}\right)^{1/q_{1}^{\prime}}\left(\fint_{Q}w_{2}^{-q_{2}^{\prime}}\right)^{1/q_{2}^{\prime}}\leq C.$$
This finishes the proof of Theorem $\ref{maximal operator for necessary condition}$. $\hfill$ $\Box$

\end{document}